\documentclass[11pt]{article}

\usepackage{hyperref}
\usepackage{mathrsfs}
\usepackage{amsfonts}
\usepackage{graphicx}
\usepackage{latexsym,bm,amsfonts,amssymb,pifont, CJK}
\usepackage[fleqn]{amsmath}
\usepackage{amsmath, color}

\usepackage{enumerate}
\newcommand{\ep}{\epsilon}

\newtheorem{theorem}{Theorem}[section]
\newtheorem{Def}[theorem]{Definition}

\newtheorem{prop}[theorem]{Proposition}

\newtheorem{lemma}[theorem]{Lemma}
\newtheorem{remark}[theorem]{Remark}

\numberwithin{equation}{section}

\def\al{\alpha}

 \def\mE{\mathbb{E}}

\def\wt{\widetilde}

\def \eref#1{\hbox{(\ref{#1})}}

\def\RR{\mathbb{R}}
\def\DD{\mathbb{D}}

\def\NN{\mathbb{N}}

\def\ZZ{\mathbb{Z}}

\def\wt{\widetilde}

\def\be{{\beta}}

\def\al{{\alpha}}

\def\be{{\beta}}

\def\ga{{\gamma}}

\def \eref#1{\hbox{(\ref{#1})}}

\begin{document}

 \title{ Crank-Nicolson scheme for stochastic differential equations driven by fractional Brownian motions }
 \date{}\;
  \author{   {\sc Yaozhong Hu$^{1,}$}\thanks{Y.  Hu is
partially supported by a grant from the Simons Foundation
\#209206.},   \;\;  {\sc Yanghui Liu$^{2}$}, \;\;  {\sc David Nualart$^{3,}$}\thanks{D. Nualart is supported by the NSF grant  DMS 1512891.
} 
 \\ 
 \\
 \small
 1. Department of Math \& Stat Sciences, University of Alberta,
Edmonton, Canada, T6G 2G1\\
\small  
2. Department of Mathematics, Purdue University, West Lafayette, Indiana, 47906 USA \,~~~~~~~\\
 \small   3. Department of Mathematics, 
	  University of Kansas, 
Lawrence, Kansas, 66045 USA~~~~~~~~~~~~ }

 \maketitle
 \begin{abstract}
We study the Crank-Nicolson scheme for stochastic differential equations (SDEs) driven by multidimensional fractional Brownian motion   $(B^{1},\dots, B^{m})$ with Hurst parameter $H\in (\frac 12,1)$. It is well-known that for    ordinary differential equations with proper conditions on the regularity of the coefficients,  the Crank-Nicolson scheme achieves a convergence rate of $n^{-2} $, regardless of the   dimension.
 In this paper  we show that,  due to   the interactions between the driving processes $ B^{1},\dots, B^{m} $,   
   the corresponding Crank-Nicolson scheme for   $m$-dimensional SDEs   has a slower   rate than for the one-dimensional SDEs. 
   Precisely, we shall prove that when  $m=1$ and  when the drift term is zero, the Crank-Nicolson scheme     achieves the exact  convergence rate $n^{-2H}$, while in the case $m=1$ and the drift term is non-zero, the exact 
   rate turns out  to be $n^{-\frac12 -H}$.
 In the general case when $m>1$, the  exact  rate equals $n^{\frac12 -2H}$.  In all these cases the limiting   distribution of the leading  error is proved to satisfy  some linear SDE
 driven by   Brownian motions   independent of the given fractional Brownian motions.
\end{abstract}

\section{Introduction}

This paper is concerned with the following   stochastic  differential equation (SDE) driven by fractional Brownian 
motion  on $ \mathbb{R}^d$
 \begin{eqnarray}
X_t  &= &x  +      \int^t_0 V  (X_s)dB_s \,, \quad
 t \in [0, T]\,,  \label{e.1.1}
 \end{eqnarray}
where
 $B = (B^{0},B^{1},\dots, B^{m})  $, and $ ( B^{1},\dots, B^{m}) $  is an $m$-dimensional  fractional Brownian motion~(fBm)   with Hurst parameter $H>\frac12$.  For notational convenience we denote  $B^{0}_{t}=t$ for $t \in [0,T]$ in order  to include the drift term  in \eref{e.1.1}. 
  The       integral on the right-hand side of \eref{e.1.1}  is of Riemann-Stieltjes type. It is well-known  that if 
 the  vector field  $   V = (V_0,V_1, \dots ,V_m): \RR^{d} \to \mathcal{L} (\RR^{m+1}, \RR^{d})  $  
   has bounded partial derivatives which are H\"older continuous of order   $  \al  >  \frac 1H - 1 $, then there exists a unique solution for equation \eref{e.1.1}, which has bounded $\frac{1}{\ga}$-variation on $[0,T]$ for any $\ga<H$; see e.g. \cite{Lyons94, Nualart}.

As in the Brownian motion case, the explicit solution of    SDEs driven by fractional Brownian motions are rarely known. Thus one has to rely on numerical methods for   simulations of these equations. 
 The simplest time-discrete  numerical approximation scheme is the  Euler scheme
 \begin{eqnarray}
X^{n}_{t_{k+1}} &=& X^{n}_{t_{k}} + V(X^{n}_{t_{k}}) (B_{t_{k+1}} - B_{t_{k}}),
\nonumber
\\
X^{n }_0 &=& x,
\label{eq1.2i}
\end{eqnarray}
where $k=0,1,\dots, n-1$ and $t_{k} = kT/n$.   This  scheme has first been  considered in \cite{Neuen, Neuenkirch} for SDEs in the one-dimensional case, and   generalized in \cite{HLN, Mishura}   to  the multidimensional case.  The solution of \eref{eq1.2i} has the exact strong convergence rate of $n^{1-2H}$ when $H>\frac12$.  When     $H=\frac12$ the scheme converges to the corresponding It\^o SDE
\begin{eqnarray*}
 X_{t} &=&x+ \int_{0}^{t}V(X_{s}) \delta B_{s} \,, \quad
 t \in [0, T],
\end{eqnarray*}
where $\delta$ denotes the It\^o  stochastic integral. 
 Note also that the Euler scheme is not convergent when $H< \frac12 $; see e.g. \cite{DNT}.   A modified Euler scheme, introduced     in \cite{HLN},   generalizes the classical Euler scheme   to the fBm case
 \begin{eqnarray}
X^{n}_{t_{k+1}} &=& X^{n}_{t_{k}} + V(X^{n}_{t_{k}}) (B_{t_{k+1}} - B_{t_{k}}) + \frac12 \sum_{ j=1}^{m} \partial V_{j}V_{j} (X^{n}_{t_{k}}) \Big(  \frac{T}{n}  \Big)^{2H} ,
\nonumber
\\
X^{n }_0 &=& x.
\label{eq1.2}
\end{eqnarray}
The modified Euler scheme has been shown to have a better convergence rate than \eref{eq1.2i}. More precisely, the rate is   $n^{\frac12 -2H}$ when $\frac12 <H<\frac34$ and  $ n^{-1} {\sqrt{\ln n}}  $ when $H = \frac34$, and in the case $\frac34 <H<1$  the rate  is  $n^{-1}$. 
  The   weak convergence rates and the asymptotic error distributions were also obtained for   this modified Euler scheme.  In \cite{GN}, the authors considered  Taylor schemes derived from the Taylor expansion  in the one-dimensional case. In \cite{HLN2}, the Taylor schemes   and their modifications were introduced    for  SDEs driven by fBms $B^{1},\dots, B^{m}$,  with Hurst parameters $H_{1},\dots, H_{m}$, where $ H_{1},\dots, H_{m} \in (\frac 12, 1] $ are not necessarily equal.  
   In \cite{D}, the Milstein scheme (or $2$nd-order Taylor scheme) has been considered for 
    the rough case $ H<\frac12$ and it is   convergent as long as $H>\frac13$.   An extension of the result to $m$-th order Taylor schemes is contained in \cite{FV}. In \cite{DNT, FR}, some $2$nd  and $3$rd order implementable schemes are studied via the Wong-Zakai approximation of~\eref{e.1.1}.  
 
  The Crank-Nicolson (or Trapezoidal) scheme has    been studied only recently. Recall that
  the Crank-Nicolson scheme for $\eref{e.1.1}$ is defined as follows:
\begin{eqnarray}
X^{n }_{t_{k+1}}  &=&  X^{n }_{t_k}+
\frac12  \big[ V   (X^n_{  t_{k+1}} ) + V   (X^n_{t_k } )\big] (B_{t_{k+1}}-B_{t_k}),\nonumber   
\\
X^{n }_0 &=& x ,
  \label{e1.2}
\end{eqnarray} 
where again $t_k=kT/n$ for $k=0,\dots, n-1$.
   In  \cite{Naganuma, Neuenkirch}, the Crank-Nicolson scheme is considered for    SDEs with   Hurst parameter $H \in (1/3, 1/2)$.   It has been shown in \cite{Neuenkirch}  that if   $V \in C^{\infty}_{b}$  the convergence rate of the Crank-Nicolson scheme   is $n^{\frac 12-3H}$.  This   rate is exact in the sense that  the renormalized error process $n^{3H-\frac12} (X-X^{n})$ converges weakly to a non-zero limit (see e.g. \cite{Naganuma}). 
   However,  due to the use of the Doss-Sussmann  representation,   these  results   are applicable  only   to the scalar SDE setting, which corresponds, with our notation, to the case 
   $m=d=1$ and $ V_{0}  \equiv 0$.   
   On the other hand, it has been conjectured in \cite{NTU} that the Crank-Nicolson scheme has exact root mean square convergence rate $n^{ \frac12-2H}$.

In view of these results, our first goal is to answer the following question:

\medskip
\noindent
{\bf Question 1:} Is the Crank-Nicolson scheme still convergent in the multidimensional setting, and is the convergence rate the same as that of  the scalar SDE? 

\medskip
Let us recall that in the case of  deterministic ordinary differential equations (ODEs), either in   the one-dimensional or multidimensional settings, and  with proper regularity assumptions on $V$, the convergence rate of the Crank-Nicolson scheme     %achieves the same convergence  rate  
is  always $ n^{-2} $. 
Surprisingly,  as we will show in this paper,    the Crank-Nicolson scheme \eref{e1.2} for equation \eref{e.1.1} has   very different features   comparing to the ODE cases.  It turns out that,
while the Crank-Nicolson scheme in the multidimensional case still converges to the solution $X$ of equation  \eref{e.1.1},  the convergence rate is largely ``throttled'' due to  the interactions between the driving processes $B^{0}$, $B^{1}$, \dots, $B^{m}$. More precisely, we will prove the following result.  Let
 $X^n_t$ denote  the continuous time   interpolation of the Crank-Nicolson scheme defined by
\begin{eqnarray}  \label{4.1}
X^{n }_t  &=&  X^{n }_{t_k}+
\frac12    \left[  V    (X^n_{t_k }  ) + V   (X^n_{  t_{k+1} }  )  \right] (B_t-B_{t_k}), 
\end{eqnarray}
for  $t \in [t_{k},t_{k+1})$, $k=0, \dots, n-1$.

\begin{theorem}\label{thm1.1}
Let $X$  be  the solution of   equation \eref{e.1.1} and  let $X^n$ be the  continuous time   interpolation of the Crank-Nicolson scheme $\left\{X^{n }_{t_0}, X^{n }_{t_1}, \cdots,  X^{n }_{t_n}\right\}$ defined by   \eref{4.1}.  Assume      that $V  \in C^{3}_{b}$.    
   Then  for any   $p \geq 1$  there exists a constant $K=K_p$ independent of $n$ such that   the following strong convergence result holds true for all~$n\in \NN$:
 \begin{eqnarray}\label{eq4.2}
\sup_{t\in [0,T]} \left( \mE |X_{t}-X^{n}_{t} |^{p } \right)^{1/p} &\leq &
K/\vartheta_{n},
\end{eqnarray}
where $\vartheta_{n}$ is defined as
\begin{eqnarray*}
\vartheta_{n}&= & 
\begin{cases} 
  n^{2H-\frac12}       & \quad \text{when}\quad m>1,
\\ 
 n^{H+\frac12}  & \quad \text{when}\quad m=1 \text{ and } V_{0} \not\equiv 0, 
\\
 n^{2H} &\quad \text{when}\quad m=1 \text{ and } V_{0} \equiv 0.
 \end{cases}
\end{eqnarray*}
\end{theorem}
Note  from Theorem \ref{thm1.1}   that if $m=1$ and $V_0\equiv 0$,  the convergence rate of the Crank-Nicolson scheme  (\ref{e1.2}) is $n^{-2H}$. 
 This result coincides with   the case of deterministic ODEs if we formally set  $H=1$, and also  the case of   SDEs driven by a one-dimensional  Brownian motion,  that is, $H=\frac12$ (see, e.g.  \cite{Naganuma, Neuenkirch}). If $m=1$ and $V_0 \not \equiv 0$, then the rate turns out to be $n^{-H-\frac 12}$. In    the general case  when $m>1$,  the converges rate becomes $n^{\frac12 -2H}$, which coincides with the modified Euler scheme defined in \eref{eq1.2} when $\frac12 <H<\frac34$. Note also that this gives a positive answer to the conjecture raised in \cite{NTU} under this general assumption.  The slowing down of convergence rate from one-dimensional case to multidimensional  cases 
% observed from Theorem \ref{thm1.1}  is captured thanks to an explicit expression of the error process $X-X^{n}$.  We will see from this explicit expression   that in the multidimensional case, the convergence of $X-X^{n}$ is dominated by 
 is due to the nonvanishing    L\'evy area term,  while in the one-dimensional case, these L\'evy area type processes disappear  and   the convergence rate of $X-X^{n}$ is then 
 dictated  by the higher order terms. 
    
    The second part of the paper is motivated by the following   question: 
    
    \medskip
    \noindent
    {\bf Question 2:} Are the convergence rates obtained in Theorem \ref{thm1.1} 
 exact? If yes,   what is 
  the  limiting  distributions   of the leading term  for  both  the one-dimensional and multidimensional cases?
  
  \medskip
 Note that the   different features observed  between the one-dimensional   and the multidimensional cases   are   true  only if the rates are exact.  
 To this aim,  we consider 
   the piecewise constant interpolations. Namely, we consider the processes  $\tilde{X}^{n}$ and $\tilde{X}$, where
 \begin{eqnarray}
 &&\tilde{X}_t :=X_{t_k}   \quad \text{and} \quad 
 \tilde{X}^n_t:=  {X}^n_{t_k} ,
 \label{5.1}
 \end{eqnarray}
 for $t\in [t_k, t_{k+1})$, $k=0,1,\dots, n$,
   and as a consequence we have
$   \tilde{X}_T:= X_T \quad \text{and} \quad \tilde{X}^n_T:= X^n_T$.   
% \begin{eqnarray}\label{e1.8}
%\tilde{X}_T:= X_T \quad \text{and} \quad \tilde{X}^n_T:= X^n_T.
%\end{eqnarray}
Recall that    here $X$ is the solution of equation \eref{e.1.1} and   $X^{n}$ is the solution of \eref{e1.2}. 
 
 The following theorem provides a complete picture of the asymptotic behaviors of the Crank-Nicolson scheme. 
    \begin{theorem}\label{thm1.2}
Let $\tilde{X}$ and $\tilde{X}^n$ be the  processes defined in (\ref{5.1}). Suppose that $V\in  C^3_b$. 
Denote $ \phi_{jj'} =   \partial  V_j V_{j'}
   -\partial  V_{j'} V_{j }$ for $j,j'=0,1,\dots,m$. 

\begin{itemize}
\item[(i)]
 Assume that      $m>1$ or $m=1$ but $V_{0} \not \equiv 0$. Then we have the convergence
 \begin{eqnarray}
( \vartheta_{n} ( \tilde{X}  -  \tilde{X}^{n}  ), B) &\rightarrow & (U,B)  \label{5.3}
\end{eqnarray}
in the Skorohod space  $D([0,T]; \mathbb{R}^{d+m+1})$ as $n$ tends to infinity. 
In the case $m>1$,    the process $U$ is the solution of the following linear SDE on $[0,T]$
\begin{eqnarray}  \label{5.4}
dU_{t}&= & \sum_{j=0}^m  
 \partial  V_j (  X_t  )U_{t} dB^{j }_t
  +
  T^{2H-\frac12}  \sqrt{  \frac{ \kappa}{2} }    \sum_{1 \leq j' < j \leq m }   \phi_{jj'} (X_{t} ) dW^{j'j}_{t}, \quad U_{0}= 0,
\end{eqnarray}
where $W=(W^{j'j})_{1\le j'<j \le m}$ is a standard $\frac {m(m-1)} 2$-dimensional Brownian motion independent of $B$ and $\kappa$ is the constant defined in \eref{e.3.1} in Section 3. In the case $m=1$ and $V_{0} \not \equiv 0$, $U$ is the solution of the  following  linear SDE on $[0,T]$
 \begin{eqnarray}
dU_{t}&= & \sum_{j=0,1}   
 \partial  V_j (  X_t  )U_{t} dB^{j }_t
  +
  T^{H+\frac12} \sqrt{  \frac{   \varrho}{2} }     \phi_{10} (X_{t} ) dW_{t}, \quad 
  U_{0}=0,  \label{5.5}
\end{eqnarray}
where $W$ is a one-dimensional standard Brownian motion independent of $B$ and $\rho$ is the constant defined in  \eref{rho} in Section 3. 
 
\item[(ii)]
Assume that    $m=1$ and  $V_{0} \equiv 0$. Then, we have the following convergence in $L^p(\Omega)$ for all $p\ge 1$ and   $t\in [0,T]$: 
\begin{eqnarray}  \label{5.6}
 n^{2H}( \tilde{X}_{t}  -  \tilde{X}^{n}_{t}  ) &\rightarrow &U_{t}\,,
 \end{eqnarray}
   where the process $U$ satisfies  the following linear SDE on $[0,T]$
\begin{eqnarray}  \label{5.7}
dU_{t}&=&    
 \partial  V  (  X_t  )U_{t} dB_t
  -\frac{T^{2H}}{4} \sum_{i,i'=1}^{d} ( V^{i}V^{i'}\partial_{i}\partial_{i'}V ) (X_{t} ) dB_{t}, \quad 
  U_{0}=0.
\end{eqnarray}
\end{itemize}
\end{theorem}

   Theorem \ref{thm1.2} shows that  in the cases   $m>1$ or $V_{0} \not\equiv 0$, one obtains  the central limit theorem   for the renormalized error process $\vartheta_{n} (X-X^{n})$, while in the case $m=0$ and $V_{0} \equiv 0$,  one gets the convergence   
    in $L^{p}$.
    %   which is associated with the situation of deterministic ODEs.
 It is interesting to point out
  that the cutoff of the convergence rates   observed in \cite{HLN, NTU} is not present in either of these cases here.      
    
  Our approach to prove    Theorem \ref{thm1.1} and Theorem \ref{thm1.2} is based on the explicit expression of $X-X^{n}$ we have mentioned previously, similar to that established in \cite{HLN}. A significant difficulty is the integrability of the Malliavin derivatives of the approximation $X^{n}$. This is due to the fact that  the Crank-Nicolson scheme \eref{e1.2} is determined by an implicit equation. 
  This difficulty will be handled        thanks to   some   fractional calculus   techniques, see e.g. \cite{CNP, HLN2, Za}.  A special attention has to   be paid also to the L\'evy area type processes mentioned  above. Our approach to handle these processes relies on a combination of fractional calculus and Malliavin calculus tools. 
   
The paper is structured as follows. In Section \ref{section 2}, we recall some basic results on the fBms as well as some upper bound estimate results and limit theorem results on fractional integrals.  
  In Section~\ref{sec3}, we consider the moment estimates and the weak convergence of  some  L\'evy area type processes. In Section~\ref{sec4}, we prove Theorem \ref{thm1.1}, and then in Section \ref{sec5}, we prove Theorem \ref{thm1.2}. Some auxiliary results are stated and proved in the appendix.

%{\color{red} 
%\begin{itemize}
%\item Properties on $\kappa$;
% \item Replace $nt/T$ by an integer?
%\item Mention $p \geq 1$;
%\item Mention $s<t$;
%\item Replace function by process;
%\item more specific reference;
%\item Specify how an application satisfies the assumption of some lemma; 
% \end{itemize}}

\section{Preliminaries}\label{section 2}

\subsection{Fractional Brownian motions}

We briefly review some basic facts about the stochastic calculus   with respect to a fBm. The reader is referred to \cite{Peccati2, N} for further details. Let $B= \{B_{t}\,, \,t\in [0,T]\}$ be a one-dimensional fBm with Hurst parameter $H \in (\frac12, 1)$, defined on some complete probability space~$(\Omega, \mathscr{F}, P)$. Namely, $B$ is a mean zero Gaussian process with covariance
\begin{eqnarray*}
\mE(B_{s}B_{t}) &=& \frac{1}{2} (t^{2H} +s^{2H}-|t-s|^{2H})
\end{eqnarray*}
for $s,t \in [0,T]$. Let $\mathcal{H}$ be the Hilbert space defined as the closure of the set of step functions on $[0,T]$ with respect to the scalar product 
\begin{eqnarray*}
\langle \mathbf{1}_{[0,t]}, \mathbf{1}_{[0,s]} \rangle_{\mathcal{H}} &=&   \frac{1}{2} (t^{2H} +s^{2H}-|t-s|^{2H}).
\end{eqnarray*}
It is easy to verify that 
\begin{eqnarray}\label{e1.4}
\langle \phi, \psi \rangle _{\mathcal{H}} &=&H(2H-1) \int_{0}^{T} \int_{0}^{T}   \phi_{u} \psi_{v} |u-v|^{2H-2} dudv
\end{eqnarray}
for every pair of step functions $ \phi, \psi \in \mathcal{H} $.

The mapping $ \mathbf{1}_{[0,t]} \mapsto B_{t} $ can be extended to a linear isometry between $\mathcal{H}$ and the Gaussian space spanned by $B$. We denote this isometry by $h \mapsto B(h)$. 
In this way, $\{ B(h), \, h \in \mathcal{H} \}$ is an isonormal Gaussian process indexed by the Hilbert space $\mathcal{H}$. 

Let $\mathcal{S}$ be the set of smooth and cylindrical random variable of the form
\begin{eqnarray*}
F&=& f(B_{t_{1}},\dots, B_{t_{N}}),
\end{eqnarray*}
where $N\geq 1$,  $t_1,\dots,t_{N} \in [0,T]$ and $f \in C^{\infty}_{b} (\RR^{N})$, namely, $f$ and all its partial derivatives are  bounded. The derivative operator $D $ on $F$ is defined as the $\mathcal{H}$-valued random variable
\begin{eqnarray*}
D_{t}F &=& \sum_{i=1}^{N} \frac{\partial f}{\partial x_{i}} ( B_{t_{1}},\dots, B_{t_{N}} ) \mathbf{1}_{[0,t_{i}]} (t), \quad t \in [0,T].
\end{eqnarray*}
For   $p\geq 1$ we define the Sobolev space $\DD^{1,p}_{B}$ (or simply $\DD^{1,p}$) as the closure of $\mathcal{S}$ with respect to the norm 
\begin{eqnarray*}
\|F\|_{\DD^{1,p}} &=& \Big( \mE[ |F|^{p}] +\mE\left[ \|DF\|_{\mathcal{H}}^{p} \right] \Big)^{1/p}.
\end{eqnarray*}
The above definition of the Sobolev space $\DD^{1,p}$ can be extended  to $\mathcal{H}$-valued random variables (see Section 1.2 in \cite{N}). We denote by $\DD^{1,p}_{B} (\mathcal{H})$ (or simply $\DD^{1,p}  (\mathcal{H})$) the corresponding Sobolev space. 

We denote by $\delta$ the adjoint of the derivative operator $D$. 
We say $u \in \text{Dom}\, \delta$ if there is a $\delta (u) \in L^{2}(\Omega)$ such that for any $F \in \DD^{1,2}$ the following duality relationship holds 
\begin{eqnarray} \label{ipf1}
\mE(\langle u, DF \rangle_{\mathcal{H}}) & =& \mE(F \delta(u)).
\end{eqnarray}
The random variable $\delta (u)$ is also called the Skorohod integral of $u$ with respect to the fBm $B$, and we use the notation $\delta (u) = \int_{0}^{T} u_{t} \delta B_{t}$\,.  The following result is an example of application of the duality relationship that will be used later in the paper.

\begin{lemma}
Let $B$ and $\wt{B}$ be independent one-dimensional fBms with Hurst parameter $H \in (\frac 12, 1)$.   Take $h \in \mathcal{H} \otimes \mathcal{H}$, then the integral $\int_{0}^{T} \int_{0}^{T} h_{s,t} \delta {B}_{s} \delta  \wt{B}_{t}$ is well defined. Denote by $D$ and $\wt{D}$ the derivative operators associated with $B$ and $\wt{B} $, respectively.  Take $F\in \DD_{\wt{B}}^{1,2}$ and assume that ${\wt{D}F \in  {\DD}_{B}^{1,2} } (\mathcal{H}) $. Then, applying the integration by parts twice, we obtain
 \begin{eqnarray} 
\label{ipf}
\mE(\langle h, D\wt{D} F \rangle_{\mathcal{H} \otimes \mathcal{H} }) &=& \mE\Big(F  \int_{0}^{T}\int_{0}^{T} h_{s,t} \delta {B}_{s} \delta  \wt{B}_{t} \Big).
\end{eqnarray}
\end{lemma}

\subsection{Weighted random sums}\label{sec2}

In this subsection, we recall some  estimates  and limit results for  Riemann-Stieltjes integrals of stochastic processes. Our main references are \cite{CNP, HLN, HLN2, Za}.  Let us start with the definition of H\"older continuous functions in  $L^p:=L^p(\Omega)$. 
In the following $\|\cdot \|_p $ denotes the $L^p$-norm in  the space $L^p$, where $p\ge 1$.

\begin{Def}
Let  $f=\{ f(t),  t \in [a,b]\}$ be a continuous   process such that $f(t) \in L^p$ for   all \,$t \in [a,b]$.  Then $f$ is called a $\be$-H\"older continuous function   in $L^p$ if the following relation holds true for all $s,t \in [ a, b ]$:
\begin{eqnarray*}
\|f(t)- f(s) \|_{p} &\le & K |t-s|^{ \beta  } .
\end{eqnarray*}
We denote by $ \| f \|_{\beta, p} $ the H\"older semi-norm
\begin{eqnarray*}
\|f\|_{  \beta, p} &=& \sup \left\{\frac{ \|f(t) - f(s)\|_p }{|t-s|^{\beta} }: t, s \in [a,b], t \not =s\right\}.
\end{eqnarray*}
 \end{Def}
  
Our first result provides an upper-bound estimate for the $L^p$-norm of a   Riemann-Stieltjes   integral. 
\begin{lemma}\label{lem 2.3}
Take
 $p\geq1$, $p',q' >1: \frac{1}{p'} +\frac{1}{q'} =1$ and $\beta,\beta' \in (0,1): \beta+\beta'>1$.
Let $f(t)$, $g(t)$, $t \in [a, b]$ be   H\"older continuous functions of order $\beta$ and $\beta'$ in $L^{pp'} $ and $L^{pq'} $, respectively.
    Then the Riemann-Stieltjes integral $ \int^b_a f  dg $ is well defined in $L^p$\,, and we have the estimate
  \begin{eqnarray}\label{eq 2.4}
\Big\|
\int^b_a f  dg  
\Big\|_p
   &\leq&
  \left( K  \|f\|_{  \beta, pp'}  
    + \|f(a)\|_{pp'} \right)  \|g\|_{  \beta', pq'} (b-a)^{\beta'},
\end{eqnarray}
where $K$ is a constant depending only on the parameters $p,p',q',\beta,\beta'$.
\end{lemma}
 \noindent \textit{Proof:}\quad
 The proof is based on  the fractional integration by parts formula (see \cite{Za}), following the arguments used in the proof of  Lemma 11.1 in \cite{HLN}. 
\hfill $\Box$

\smallskip

 Given a double sequence of random variables $\zeta= \{ \zeta_{k, n} , n \in \mathbb{N}, k=0, 1, \dots, n \}$, for each $t \in [0,T]$ we set
\begin{eqnarray}\label{e7.1}
g_n(t): = \sum_{k=0}^{   \left\lfloor \frac{nt}{T} \right\rfloor } \zeta_{k , n} \,,
\end{eqnarray}
where $\left\lfloor \frac{nt}{T} \right\rfloor$ denotes the integer part of $\frac {nt}T$. 
 	 We recall the following result   from    \cite{HLN2}, which
  provides an upper-bound estimate for  weighted random sums (or the so-called discrete integrals) of the process $g_{n}$.
   \begin{lemma}\label{lem 2.4}
 Let $p$, $p'$, $q'$, $\beta$, $\beta'$ be as in Lemma \ref{lem 2.3}.  Let $f$ be  a H\"older continuous function of order $\beta$ in $L^{pp'}$. Let $g_n $  be as in \eref{e7.1} such that for any $j,k =0, 1, \dots, n $ we have
\begin{eqnarray*} 
 \mE (| g_n( kT/n ) - g_n ( jT/n)  |^{pq'} ) &\leq&
  K ( |k-j|/n)^{ \beta' pq'}.
 \end{eqnarray*}
   Then the following estimate holds true for $i,j=0, 1,\dots, n $, $ i>j$:
\begin{eqnarray*}
\Big\| \sum_{k= j  +1   }^{ i   } f(t_k) \zeta_{k, n} \Big\|_p &\leq& K \left(     \|f\|_{\beta, pp'}   +   \|f(t_{j} )\|_{pp'} \right)
 \Big(\frac{i-j}{n}\Big)^{\beta'} .
\end{eqnarray*}
\end{lemma}
   
\medskip

Let us now recall    some limit theorems  for weighted random sums. The first result says that if the ``weight-free'' random sum \eref{e7.1} converges weakly and if the weight process satisfies certain regularity  assumption, then the weighted random sum also converges weakly. The reader is referred to  \cite{CNP} for further details. 

\begin{prop}\label{prop2.8}
Let   $g_n$ be   defined in \eref{e7.1}. Assume   that $g_{n}$ satisfies the inequality
\begin{eqnarray*} 
 \mE (| g_n( kT/n ) - g_n ( jT/n)  |^{4} ) &\leq&
  K ( |k-j|/n)^{ 2}
 \end{eqnarray*}
 for  $j,k =0, 1, \dots, n $.
 Suppose further that the finite dimensional distributions of 
$g_n$ converge stably to those of $W=\{W_{t}, t\in [0,T]\}$, where $W$ is a standard Brownian motion   independent of $g_{n}$.    Let $f=\{f(t), t\in [0,T]\} $ be a $\beta$-H\"older continuous process for $\beta>1/2$. 
Then the finite dimensional distributions of 
$\sum_{k=0}^{\left\lfloor \frac{nt}{T} \right\rfloor } f(t_{k}) \zeta_{k, n} 
$ converge stably to those of 
$\int_{0}^{t} f(s) dW_{s}$, { where recall that $\left\lfloor \frac{nt}{T} \right\rfloor$ denotes the integer part of $\frac {nt}T$.}
 \end{prop}
Recall that a sequence of random vectors $F_{n}$ converges stably  to a random vector $F$, where $F$ is defined on an extension $(\Omega', \mathscr{F}', \mathbb{P}')$ of the original probability $(\Omega, \mathscr{F}, \mathbb{P})$, if $(F_{n}, Z) \rightarrow (F,Z)$ weakly for any $\mathscr{F}$-measurable random variable $Z$. The reader is  referred to \cite{Aldous, Jacod, Renyi} for further details on stable convergence.

\medskip

The following result  can be viewed as the $L^p$-convergence version of Proposition \ref{prop2.8}; see \cite{HLN}.

 \begin{prop}\label{prop7.1}
 Take $\beta, \lambda \in (0,1):  \beta +\lambda>1 $.
   Let $p\geq1$ and $p', q' > 1$ such that $\frac{1}{p'} +\frac{1}{q'} =1$ and $pp'>\frac{1}{ \beta}$, $pq' > \frac{1}{\lambda}$. Let $g_n $  be defined in \eref{e7.1}. Suppose that the following two conditions hold true:
\begin{itemize}
\item[(i)] For   $t \in [0, T]$, we have the convergence $g_n(t ) \rightarrow z(t)$ in $L^{pq'}$\,;
\item[(ii)]     For   $j,k =0, 1, \dots, n $ we have the relation:
\begin{eqnarray*}
 \mE (| g_n( kT/n ) - g_n ( jT/n)  |^{pq'} ) & \leq& K ( |k-j|/n)^{\lambda pq' }.
 \end{eqnarray*}
\end{itemize}
Let $f=\{ f(t) , t \in [0, T]\}$  be  a  continuous process  such that $\mE (\|f\|_{\beta}^{pp'} ) \leq K$ and $\mE (|f( 0 )|^{pp'}) \leq K$.
 Then for each $t\in [0, T]$ we have the convergence:
\begin{eqnarray*}%\label{eqn7.3}
 \lim_{n\to \infty}  \sum_{k=0}^{   \left\lfloor \frac {nt}T \right\rfloor }f(t_k) \zeta_{k , n} = \int_0^t f(s) dz(s)
   \,, 
\end{eqnarray*}
where the limit is  understand as the limit in $L^p$. 
\end{prop}

\section{L\'evy area type processes}\label{sec3}

Let   $B= \{ B_{t}, t \geq 0 \}$ be a  one-dimensional fBm with Hurst parameter $H \in (\frac12, 1)$, and let   $\wt{B} = \{ \wt{B}_{t} , t \geq 0  \}$ be a H\"older continuous  process of order $\beta>\frac12$.  Let  $\Pi= \{0=t_{0}<t_{1}<\cdots<t_{n}=T \}$ be the uniform partition on $[0,T]$ and take $t_{n+1} =   \frac {n+1}{n}T$.
 We consider the following L\'evy area type process on $[0,T]$
 \begin{eqnarray}\label{e 3.1}
  {Z}_n(t)&=&
  \sum_{k=0}^{ l  } \left(
    \int_{t_k}^{t_{k+1}   } \int_{t_k}^{s } d\wt{B}_u  dB_s -
   \int_{t_k}^{t_{k+1}    } \int^{t_{k+1}   }_{s } d\wt{B}_u dB_s \right)
 \,
 \end{eqnarray}
 for $t\in [t_{l}, t_{l+1})\cap [0,T]$, $l=0,\dots, n $. 
  In this section, we study the convergence rate and the asymptotic distribution of the sequence $\{Z_{n},  n \in \NN\}$. We will mainly focus   on       two cases: 
  (i)   $\wt{B}$ is  an    independent copy of $B$; and (ii) $\wt{B}$ is the identity function, that is, $\wt{B}_t = t $ for $t \geq 0$.

\subsection{Case (i)}
 For simplicity, we denote by  $\mu$ the measure on the plane $\mathbb{R}^2 $ given by 
 \begin{eqnarray*}
\mu (dsdt) = H(2H-1)|s-t|^{2H-2}dsdt.
\end{eqnarray*}
  For each $p\in \ZZ$ we set  
  \begin{eqnarray*}
      Q( p)
       =          \int_{ {0}}^{ { 1}}  \int_{ { p} }^{ {{ p} +1}}
      \int_{ {0}}^{  t }  \int_{ { p} }^{s}
       \mu(dvdu) \mu(ds dt)
         ,
             \quad
        R(p) =
          \int_{ {0}}^{ { 1}}  \int_{ p}^{ { p+1}}
      \int_{ {t}}^{  1 }  \int_{ p}^{s}
       \mu(dv du) \mu(ds dt).
    \end{eqnarray*}
The following   result provides some properties of the process $Z_{n}$.
 \begin{prop}\label{prop 3.1}
 Let $Z_n$ be the process defined in \eref{e 3.1} and let $\wt{B}$ be an independent copy of $B$. Then, there exists a constant $K$ depending on $H$ and $T$ such that for $t,s \in \Pi $ we have
\begin{eqnarray}\label{e 3.1i}
 n^{4H-1} \mE( [  {Z}_n(t)-  {Z}_n(s) ]^2 ) &\leq& K |t-s|.
 \end{eqnarray}
Furthermore,
%the finite dimensional distribution of  $n^{2H-\frac12} Z_n$ converges stably to that of $  \sqrt{  \frac{2\kappa}{T} } W$ as $n$ tends to infinity, in other words, 
the finite dimensional distributions of 
$
\left( n^{2H-\frac12} Z_n(t)  , B_{t} , \wt{B}_{t}  , \, t\in [0,T]  \right) 
$
converge weakly to those of
$
\left(   T^{2H-\frac12} \sqrt{  {2\kappa}  }   W_{t}  , B_{t}  , \wt{B}_{t}, \, t\in [0,T]  \right)
 $
   as $n$ tends to infinity,
where $W= \{W_t, t \in [0, T]\}$ is a standard Brownian motion independent of $(B,\wt{B})$,  and
\begin{eqnarray}\label{e.3.1}
  \kappa  &=& \sum_{p \in \,\mathbb{Z} }  ( Q( p)- R(p) )   .
 \end{eqnarray}
   \end{prop}

\begin{remark}
 The following figure provides the graph of the parameter $\kappa$ as a function of $H$ on $(\frac12 , 1)$.  We observe  that $\kappa$ converges to $\frac12$ as $H$ tends to $\frac12$ which corresponds to the Brownian motion, and  $\kappa$ approaches $0$ when $H$ tends to $1$.

\includegraphics[width=10cm]{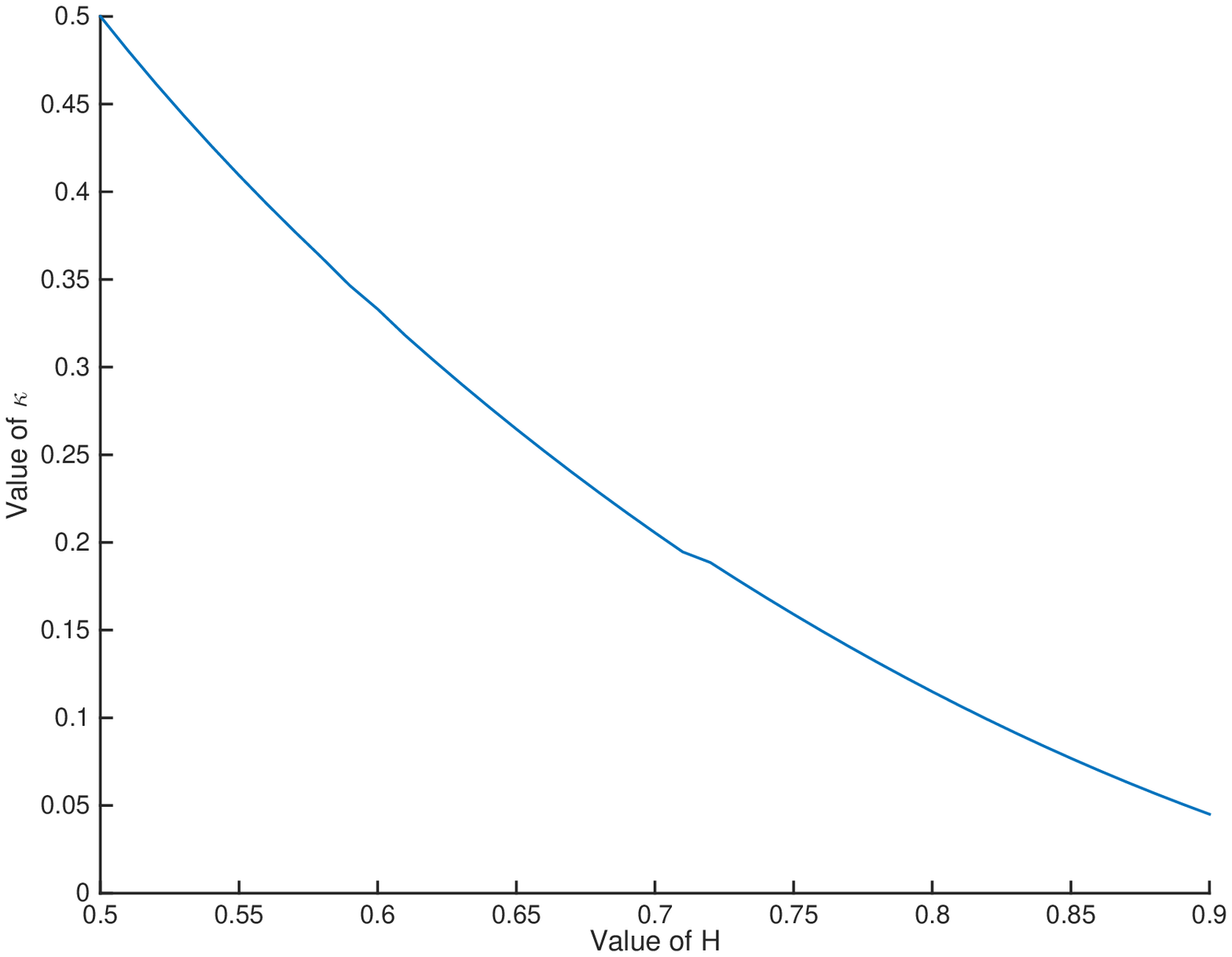}
\end{remark}

 \noindent \textit{Proof of Proposition \ref{prop 3.1}:}\quad
 The proof is divided into four  steps. 
 
 \noindent \textit{Step 1.}\quad
 In this step, we show the convergence  of   $n^{ 4H - 1 } \mE ({ {Z}}_n(t)^2 )$\,  and derive its limit as $n\rightarrow \infty$. We first calculate the second moment of $Z_{n}(t)$.    Note that when $\wt{B}$ is an independent copy of $B$ we have
 \begin{eqnarray}\label{e3.4i}
 Z_{n}(t) &=&
     \sum_{k=0}^{ \lfloor \frac{nt}{T} \rfloor  } \left(
    \int_{t_k}^{t_{k+1}   } \int_{t_k}^{s } \delta \wt{B}_u  \delta B_s -
   \int_{t_k}^{t_{k+1}    } \int^{t_{k+1}   }_{s } \delta \wt{B}_u \delta B_s \right)
\nonumber
   \\
  &=&    \sum_{k=0}^{ \lfloor \frac{nt}{T} \rfloor  } \int_{0}^{T} \int_{0}^{T}      \beta_{\frac kn} (s)    \gamma_{t_k, s} (u)       \delta \wt{B}_{u} \delta {B}_{s}\,,
\end{eqnarray}
 where $\delta$ denotes the Skorohod integral and     
 \begin{eqnarray}\label{e3.4}
 \beta_{\frac kn} (s) = \mathbf{1}_{[t_k, t_{k+1}] } (s), \quad\quad\quad\quad \gamma_{t_k, s} (u) =
    \mathbf{1}_{[t_k, s] }(u)  - \mathbf{1}_{[s, t_{k+1}] }(u).
\end{eqnarray}
   By  the integration by parts formula (\ref{ipf}) and taking into account the expression of $Z_{n}(t)$ in \eref{e3.4i} we obtain
   \begin{eqnarray}\label{e3.6}
\mE
[{ {Z}}_n(t)^2 ]
&=&  \sum_{k =0}^{ \lfloor \frac{nt}{T} \rfloor  } \int_{0}^{T} \int_{0}^{T} \int_{0}^{T}\int_{0}^{T} \wt{D}_{u'}D_{s'} Z_{n}(t)   \beta_{\frac kn} (s)   \gamma_{t_k, s} (u)   \mu(dudu') \mu(dsds')
,
\end{eqnarray}
where $D$ and $\wt{D}$ are the  derivative operators associated with $B$ and $\wt{B}$, respectively. 
It is clear that 
\begin{eqnarray*}
\wt{D}_{u'}D_{s'} Z_{n}(t) &=& 
\sum_{k=0}^{ \lfloor \frac{nt}{T} \rfloor  }        \beta_{\frac kn} (s')    \gamma_{t_k, s'} (u') .
\end{eqnarray*}
Therefore, we obtain the expression
   \begin{eqnarray*}
\mE
[{ {Z}}_n(t)^2 ]
&=& 
\sum_{k,k'=0}^{ \lfloor \frac{nt}{T} \rfloor  } \int_{0}^{T} \int_{0}^{T} \int_{0}^{T}\int_{0}^{T}  \beta_{\frac {k'}n} (s') \beta_{\frac kn} (s) \gamma_{t_{k'}, s'} (u') \gamma_{t_k, s} (u)   \mu(dudu') \mu(dsds')
.
\end{eqnarray*}
By changing the variables from $(u,u',s,s')$ to $\frac{T}{n} (u,u',s,s') $ we obtain 
   \begin{eqnarray*}
\mE[{ {Z}}_n(t)^2 ]
&=& \left( \frac{T}{n} \right)^{4H} \sum_{k,k'=0}^{ \lfloor \frac{nt}{T} \rfloor  }  
\int_{ {k'}}^{ {k'+1}} \int_{ {k}}^{ {k+1}} \int_{0}^{n}\int_{0}^{n}    \varphi_{ {k'}, s'} (u') \varphi_{ k, s} (u)   \mu(dudu') \mu(dsds'),
\end{eqnarray*}
where
$
\varphi_{ {k},s}(u) = \mathbf{1}_{[k,s]}(u) - \mathbf{1}_{[s,k+1]}(u).
$
  Denote $\varphi_{ {k},s}^{0}(u) =  \mathbf{1}_{[k,s]}(u) $, $ \varphi_{ {k},s}^{1}(u) = \mathbf{1}_{[s,k+1]}(u) $, and set
  \begin{eqnarray*}
e_{ij}&=&\int_{ {k'}}^{ {k'+1}} \int_{ {k}}^{ {k+1}} \int_{0}^{n}\int_{0}^{n}    \varphi_{ {k'}, s'}^{i} (u') \varphi_{ k, s}^{j} (u)   \mu(dudu') \mu(dsds').
\end{eqnarray*}
   Then we can write
    \begin{eqnarray*}
\mE[{ {Z}}_n(t)^2 ]
&=&  \left( \frac{T}{n} \right)^{4H}  \sum_{k,k'=0}^{ \lfloor \frac{nt}{T} \rfloor  }  \sum_{i,j=0,1}  (-1)^{i+j}e_{ij}.
\end{eqnarray*}
  It is easy to see that $e_{00}=e_{11}=Q(k-k')$ and $e_{01}=e_{10}=R(k-k')$.
 Therefore,
  \begin{eqnarray}\label{e 2.4i} 
\mE
[{ {Z}}_n(t)^2 ]
&=&
2   \left( \frac{T}{n}\right)^{4H} \sum_{k,k' =0}^{ \lfloor \frac{nt}{T} \rfloor  } [Q(k-k')-R(k-k')]
 \,.
   \end{eqnarray}
Taking $p =k-k'$ on the right-hand side of \eref{e 2.4i},    we obtain  
    \begin{eqnarray}\label{e 2.4}
\mE
[{ {Z}}_n(t)^2 ]
&= &
  2    \left( \frac{T}{n}\right)^{4H} 
  \left(
 \sum_{p =0}^{ \lfloor \frac{nt}{T} \rfloor   }\sum_{ k =p}^{\lfloor \frac{nt}{T} \rfloor  }
[Q(p)-R(p)]
 +
 \sum^{-1}_{p= - \lfloor \frac{nt}{T} \rfloor  }\sum_{ k =0}^{\lfloor \frac{nt}{T} \rfloor  +p}
[Q(p)-R(p)]
\right)
\nonumber
\\
&:= &
q_1+q_2
\,.
   \end{eqnarray}
        We decompose $q_{1} $ as follows,
    \begin{eqnarray*}
q_1 &= &
     2    \left( \frac{T}{n}\right)^{4H} 
       \sum_{p=0 }^{\lfloor \frac{nt}{T} \rfloor   }
  (\lfloor \frac{nt}{T} \rfloor  -p+1 )
  (Q(p) - R(p) )
    \\
    &= &
    2    \left( \frac{T}{n}\right)^{4H} 
\left(
\lfloor \frac{nt}{T} \rfloor  
  \sum_{p=0 }^{\lfloor \frac{nt}{T} \rfloor }
          (Q(p) - R(p) )
    -
     \sum_{p=0 }^{\lfloor \frac{nt}{T} \rfloor }
(  p-1)
    (Q(p) - R(p) )
\right)
\\
&:= &
q_{11}+q_{12}\,.
   \end{eqnarray*}
 By mean value theorem, it is easy to show that $|Q(p) - R(p) |  \leq  Kp^{4H-5} $ for $p>0$, so the sum
$    \sum_{p=0}^{\infty} ( Q(p)- R(p) )  $
is convergent, and
  that
 \begin{eqnarray*}
  \left|
    \sum_{p=0 }^{\lfloor \frac{nt}{T} \rfloor }
  (p-1)
    (Q(p) - R(p) ) \right|
    & \leq & K(n^{4H-3} \vee 1)  .
\end{eqnarray*}
   Here $a\vee b$ denotes the maximum of $a$ and $b$.    Therefore,
         \begin{eqnarray}\label{e 2.6i}
    \lim_{n \rightarrow \infty} n^{4H-1} q_{11}
    &=&
     \lim_{n \rightarrow \infty} 2n^{4H-1}       \left( \frac{T}{n}\right)^{4H}  \lfloor \frac{nt}{T} \rfloor
  \sum_{p=0 }^{\lfloor \frac{nt}{T} \rfloor }
         (Q(p) - R(p) )
         \nonumber
         \\
         &=&
        2 t T^{4H-1}    \sum_{p=0 }^{\infty}
         (Q(p) - R(p) )  
    \end{eqnarray}
and
    \begin{eqnarray}\label{e 2.6ii}
    \lim_{n \rightarrow \infty} n^{4H-1} q_{12}
    &=  &
    0 .
    \end{eqnarray}
    In summary, from \eref{e 2.6i} and \eref{e 2.6ii}, we obtain
    \begin{eqnarray}\label{e 2.6}
    \lim_{n \rightarrow \infty} n^{4H-1} q_1
  &  =  &
      2 t  T^{4H-1}   \sum_{p=0 }^{\infty}
         (Q(p) - R(p) )  .
    \end{eqnarray}
  In a similar way,   we can prove the following convergence for $q_{2}$ 
     \begin{eqnarray}\label{e 2.7}
    \lim_{n\rightarrow \infty} n^{ 4H - 1} q_2 &=&   2 t T^{4H-1}  \sum^{-1}_{p= -\infty} ( Q( p)- R(p) )  .
    \end{eqnarray}
    Substituting \eref{e 2.6} and  \eref{e 2.7}
 into \eref{e 2.4} yields
 \begin{eqnarray}\label{e 2.8}
  \lim_{ n \rightarrow \infty} n^{ 4H - 1 } \mE ({ {Z}}_n(t)^2 ) &=&    2     T^{4H-1} \kappa t,
  \end{eqnarray}
where recall that $\kappa$ is a constant defined in \eref{e.3.1}.

\medskip

\noindent \textit{Step 2.}\quad 
In this step, we show the inequality \eref{e 3.1i}.   This inequality  is obvious when $s=t$. In the following, we consider the case when $t>s$.

Take $t\in \Pi$. By the definition of $q_{1}$ we have
 \begin{eqnarray*}
q_1 &\leq &
     2    \left( \frac{T}{n}\right)^{4H} 
       \sum_{p=0 }^{  \frac{nt }{T}   }
  (   \frac{nt}{T}    -p+1)
\,\,  | Q(p) - R(p) |
  \\
  &\leq&
   2     \left( \frac{T}{n}\right)^{4H} 
     ( \frac{nt}{T}  +1)
       \sum_{p=0 }^{ \infty }
  | Q(p) - R(p) |\,.
      \end{eqnarray*}
  In the same way, we can show that
 \begin{eqnarray*}
q_2  & \leq&
   2     \left( \frac{T}{n}\right)^{4H}      ( \frac{nt}{T}  +1) 
       \sum_{p= - \infty }^{ -1 }
  | Q(p) - R(p) |.
      \end{eqnarray*}
Applying these two inequalities to   \eref{e 2.4} we obtain
\begin{eqnarray}\label{e3.13}
 n^{4H-1} \mE({  {Z}}_n(t)^2 ) &\leq& K (t+\frac{T}{n})
 \end{eqnarray}
for $t\in \Pi$, where $K$ is a constant depending on $H$, $T$.    
 Take $s, t \in \Pi: s<t$.  The inequality \eref{e 3.1i} then follows by replacing $t$ in \eref{e3.13} by $t-s-\frac{T}{n}$ and noticing that $Z_{n} (t) -Z_{n}(s)$ and $Z_{n} (t-s-\frac{T}{n})$ are equal in distribution and thus have the same second moments.

\medskip

\noindent \textit{Step 3.}\quad  Take $s,t \in [0,T]: s<t$.  In this step, we derive the limit of the quantity $n^{4H-1} \mE (Z_{n}(t) Z_{n}(s))$. 
Denote $\eta(t) =t_{k}$ for $t \in [t_{k},t_{k+1})$, $k=0,1,\dots, n$. Then we have $Z_{n}(t) = Z_{n}(\eta(t))$. 
Since $Z_{n} (\eta(t)) - Z_{n}(\eta(s))$ and $Z_{n} (\eta(t)-\eta(s)-\frac{T}{n}) $ have the same distribution, we have
  \begin{eqnarray}\label{e3.12i}
   \mE (|Z_{n} (t) - Z_{n}(s) |^{2}) &=&  \mE (|Z_{n} (\eta(t)) - Z_{n}(\eta(s)) |^{2})
   \nonumber
\\
 &=&    \mE (|Z_{n} (\eta(t)-\eta(s) -\frac{T}{n}) |^{2}).
\end{eqnarray}  
   Note that 
  % \begin{align*}
$0< (t-s)  - (\eta(t)-\eta(s) - \frac{T}{n})    < \, 2\frac{T}{n}$,
%\end{align*}
 so either $Z_{n} (\eta(t)-\eta(s) - \frac{T}{n}) = Z_{n}(t-s)$ or $Z_{n} (\eta(t)-\eta(s) - \frac{T}{n}) = Z_{n}(t-s - \frac{T}{n})$.
 In both cases we have
  \begin{eqnarray}\label{e3.11i}
\lim_{n\rightarrow \infty} n^{4H-1 } \mE (|Z_{n} (\eta(t) -\eta(s)  - \frac{T}{n} ) |^{2}) 
 &=& \lim_{n\rightarrow \infty} n^{4H-1 } \mE (|Z_{n} (t -s  ) |^{2}).
%\nonumber \\
% =& 2\frac{t-s}{T} \kappa. 
\end{eqnarray}
 Indeed, the identity is clear in the first case. In the second case, it can be shown with the help of  \eref{e3.13} that
 \begin{eqnarray*}
 \lim_{n\rightarrow \infty} n^{4H-1 } \left( \mE (|Z_{n} (t -s -\frac{T}{n}  ) |^{2}) -  \mE (|Z_{n} (t -s  ) |^{2})  \right) &=&0.
\end{eqnarray*}
 The identity \eref{e3.11i} then follows.
 Substituting \eref{e3.11i} into \eref{e3.12i} and with the help of \eref{e 2.8}  we obtain
  \begin{eqnarray}\label{e3.11ii}
\lim_{n\rightarrow \infty} n^{4H-1 } \mE (|Z_{n} (t) - Z_{n} (s)  |^{2}) 
 &=& 2 {T^{4H-1}} \kappa (t-s). 
\end{eqnarray}
  By expanding the left-hand side of \eref{e3.11ii} and using \eref{e 2.8}, we obtain 
  \begin{eqnarray}\label{e 3.8i}
  \lim_{ n \rightarrow \infty} n^{4H - 1} \mE [  {Z}_n(t)    {Z}_n(s)   ] 
  &=   2 {T^{4H-1}} \kappa     {(t\wedge s)} , \quad \,  s,t \in [0,T].
  \end{eqnarray}

\medskip

\noindent \textit{Step 4.}\quad 
  In this step, we    prove   the weak convergence for the finite dimensional distributions of $ (n^{2H-\frac12}Z_{n}, B, \wt{B}) $. Given $r_1, \dots, r_L \in [0, T]$, $L \in \NN$, we need  to show that the random vector
\begin{eqnarray*}
\Theta^n_L
:=\left( n^{2H-\frac12} (Z_n(r_1), \dots, Z_n(r_L)), B_{r_1}, \dots, B_{r_L}, \wt{B}_{r_1}, \dots, \wt{B}_{r_L} \right)
\end{eqnarray*}
  converges in law to
 \begin{eqnarray*}
\Theta_L
:=  \left( T^{2H-\frac12}\sqrt{  {2\kappa}  }   \left(   W(r_1), \dots,  W(r_L) \right), B_{r_1}, \dots, B_{r_L}, \wt{B}_{r_1}, \dots, \wt{B}_{r_L} \right)
\end{eqnarray*}
  as $n $ tends to infinity, where recall that   $W= \{W_t, t \in [0, T]\}$ is a standard Brownian motion independent of $(B,\wt{B})$.
     According to  \cite{Peccati} (see also Theorem 6.2.3 in \cite{Peccati2}), this is true if we can  show the weak convergence
     of  each component of $\Theta^n_L$ to  the corresponding component of $\Theta_L $
      and the convergence of its  covariance matrix to that of $\Theta_L $.

 The convergence of the covariance of $n^{2H-\frac{1}{2}  }Z_n (r_i) $ and $ n^{2H-\frac{1}{2} }Z_n (r_j) $   to that of $   T^{2H-\frac12}\sqrt{  {2\kappa}  }    W(r_i ) $ and $    T^{2H-\frac12}\sqrt{  {2\kappa}  }    W(r_j ) $ follows from \eref{e 3.8i}. The covariance of $n^{2H-\frac{1}{2} }Z_{n}(r_{i}) $ and $(B_{r_{j}}, \wt{B}_{r_{j}})$ is zero since they are in different chaos, so the limit of the covariance is zero, which equals   the covariance of $T^{2H-\frac12}\sqrt{  {2\kappa}  }   W(r_{i})$ and $(B_{r_{j}}, \wt{B}_{r_{j}})$ since $W$ and $B$ are independent.

   By the fourth moment theorem (see \cite{Nualart3} and also Theorem 5.2.7 in \cite{Peccati2}) and taking into account~{\eref{e 3.8i}}, to show the weak convergence of the components of $\Theta^n_L$ it remains to  show that the limits of their  fourth moments exist,  and 
   \begin{eqnarray}\label{e3.11}
\lim_{n\rightarrow \infty} n^{8H-2} \mE\left[     Z_{n} (t)^{4} \right] 
&=&
 3 \lim_{n\rightarrow \infty} n^{8H-2} \left(    \mE  [    Z_{n} (t)^{2}  ]  \right)^{2}
\end{eqnarray} 
for $t\in [0,T]$.

Applying    the  integration by parts   formula (\ref{ipf})   to $\mE
[{ {Z}}_n(t)^4 ]$ and taking into account the expression of $Z_{n}(t)$ in \eref{e3.4i}, we obtain
   \begin{eqnarray}\label{e.3.11}
\mE
[{ {Z}}_n(t)^4 ]
&=& \mE
[{ {Z}}_n(t)^{3} \cdot { {Z}}_n(t)  ]
\nonumber
\\
&=&
\sum_{k=0}^{ \lfloor \frac{nt}{T} \rfloor }
\mE
 \int_0^T\int_0^T\int_0^T\int_0^T
  \wt{D}_{u'}D_{s'}
\left[ { {Z}}_n(t)^3 \right]
  \beta_{\frac kn }(s) \gamma_{t_k, s} (u)
 \mu(dudu')  \mu(dsds')
  ,
   \end{eqnarray}
   where $D$ and $\wt{D}$ are the differential operators associated with   $B$ and $\wt{B}$, respectively. 
  We expand the second derivative $ \wt{D}_{u'}D_{s'}
\left[ { {Z}}_n(t)^3 \right] $ as follows
\begin{eqnarray*}
\wt{D}_{u'}D_{s'}
\left[ { {Z}}_n(t)^3 \right] &=  &3 Z_{n}(t)^{2} \wt{D}_{u'}D_{s'} Z_{n}(t)+ 6  Z_{n}(t)  \wt{D}_{u'}Z_{n}(t) D_{s'} Z_{n}(t).
\end{eqnarray*}
Substituting the above identity into \eref{e.3.11}, we obtain 
\begin{eqnarray*}
\mE
[{ {Z}}_n(t)^4 ] &= & d_{1} +d_{2}\,,
\end{eqnarray*}
where
 \begin{eqnarray}\label{eqn 3.8}
 d_2&=& 6   \sum_{k=0}^{ \lfloor \frac{nt}{T} \rfloor }
\mE
 \int_0^T\int_{0}^{T}\int_0^T\int_0^T
  { {Z}}_n(t)
  \wt{D}_{u'}
   { {Z}}_n(t)
  D_{s'}
  { {Z}}_n(t)
  \beta_{\frac kn }(s) 
    \gamma_{t_k, s} (u)
   \mu(dudu')  \mu(dsds')
   \end{eqnarray}
   and
  \begin{eqnarray} \label{e3.20} 
d_1&=&
 3 \mE [ { {Z}}_n(t) ^2 ] \sum_{k=0}^{ \lfloor \frac{nt}{T} \rfloor }
 \int_0^T \int_{0}^{T} \int_0^T\int_0^T
  \wt{D}_{u'}D_{s'}
 { {Z}}_n(t)
 \beta_{\frac kn }(s)   \gamma_{t_k, s} (u)
   \mu(dudu')  \mu(dsds').  
  \end{eqnarray}
     Substituting \eref{e3.6} into $d_{1}$, we obtain  
  \begin{eqnarray}\label{e3.14i} 
d_{1} &=&   3 
\mE [ { {Z}}_n(t) ^2 ]
\mE [ { {Z}}_n(t) ^2 ]. 
\end{eqnarray}
The term $d_2$ is more sophisticated. We shall prove in  
  Section \ref{section 6.1}  the following fact: 
\begin{eqnarray}\label{e3.14}
\lim_{n\rightarrow \infty} n^{8H-2}d_2 &=& 0.
\end{eqnarray}
 The identity \eref{e3.14i} and the convergence \eref{e3.14} together  imply the identity \eref{e3.11}.   This completes the proof.
\hfill $\Box$

\subsection{Case (ii)}

In this subsection, we consider the process $Z_n$   in \eref{e 3.1} under the assumption that    $\wt{B}_t = t$, $t\in [0, T]$. We   denote $z_n:=Z_n$ in this subsection to distinguish it from  the  $Z_{n}$ in the   previous subsection. 
For each $p \in \mathbb{Z} $, we denote
  \begin{eqnarray*}
     \wt{Q}( p)
       = 
         \int_{ {0}}^{ { 1}}  \int_{ { p} }^{ {{ p} +1}}
      \int_{ {0}}^{  t }  \int_{ { p} }^{s}
        dvdu  \mu(ds dt)
         ,
             ~~~~
        \wt{R}(p) =
          \int_{ {0}}^{ { 1}}  \int_{ p}^{ { p+1}}
      \int_{ {t}}^{  1 }  \int_{ p}^{s}
        dv du  \mu(ds dt),
    \end{eqnarray*}
  where
    recall that  $\mu (dsdt)=H(2H-1)|s-t|^{2H-2}dsdt$ is a measure on $\RR^{2}$.
\begin{prop}\label{prop 3.2}
 Let $z_n$ be the process defined in \eref{e 3.1} with $\wt{B}_{t} =t$, $t\in [0,T]$. Then, there exists a constant $K$ depending on $H$ and $T$ such that for $t,s \in \Pi = \{\frac{T}{n}i: i=0,\dots, n\}$ we have
\begin{eqnarray}\label{eq 3.17}
 n^{2H+1} \mE( [  {z}_n(t)-  {z}_n(s) ]^2 ) & \leq & K |t-s|.
 \end{eqnarray}
Moreover,
   the finite dimensional distributions of 
the process $ ( n^{ H + \frac12} z_n \,, \,B   ) 
$
converge weakly to those of  
$
  (   {  \sqrt{2  \varrho}{ } }T^{H+\frac12} W \,, \,B   )$
 as $n \rightarrow \infty$, where $W= \{W_t, t \in [0, T]\}$ is a standard Brownian motion independent of $B$, and
  \begin{eqnarray}
 \varrho  &:= & \sum_{p \in \mathbb{Z} } ( \wt{Q}(p) - \wt{R}(p) )   .  \label{rho}
\end{eqnarray}
  \end{prop}

\begin{remark}
 The following figure provides the graph of the parameter $\rho$ as a function of $H$ on $(\frac12 , 1)$.  We see  that $\rho$ converges to $\frac16$ as $H$ tends to $\frac12$, and  $\rho$ approaches $0$ when $H$ tends to $1$.

\includegraphics[width=10cm]{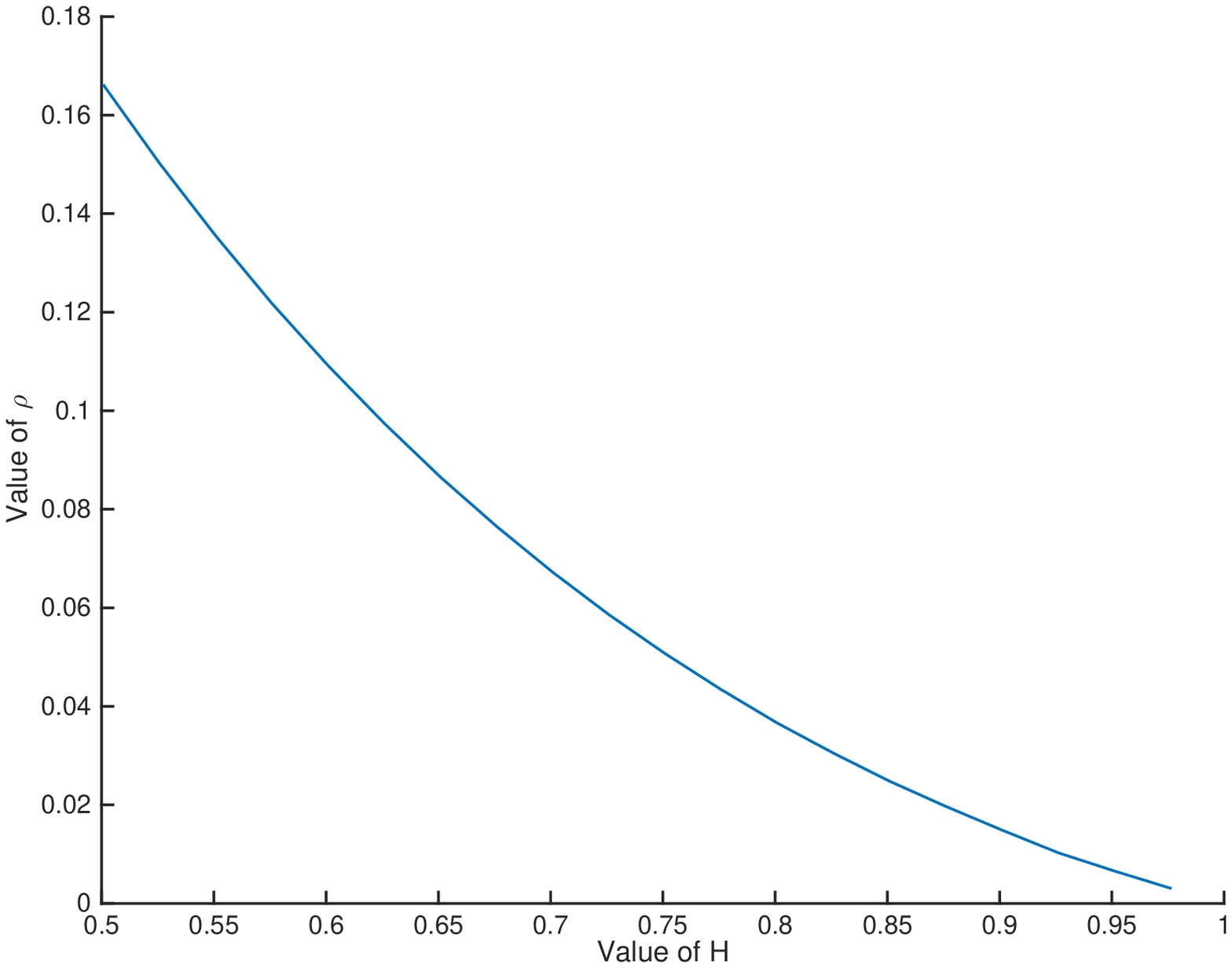}
\end{remark}

 \noindent \textit{Proof of Proposition \ref{prop 3.2}:}
 \quad 
 The proof is done in  three  steps. 
 
  \noindent \textit{Step 1.}\quad
 In this step, we   calculate the second moment of $z_{n}(t)$. Write $z_{n}(t)$ as
 \begin{eqnarray*}
z_{n}(t)&=& 
   \sum_{k=0}^{ \lfloor \frac{nt}{T} \rfloor  } \int_{0}^{T} \int_{0}^{T}      \beta_{\frac kn} (s)    \gamma_{t_k, s} (u)     d{u} \delta {B}_{s}\,,
\end{eqnarray*}
where $  \beta_{\frac kn} (s)   $ and $ \gamma_{t_k, s} (u) $ are defined  in \eref{e3.4}.
  Then, applying the covariance formula  \eref{e1.4}, we end up with
\begin{eqnarray*}
 \mE
[{ {z}}_n(t)^2 ]
 &=&  \sum_{k,k' =0}^{ \lfloor \frac{nt}{T} \rfloor } \int_{0}^{T} \int_{0}^{T} \int_{0}^{T} \int_{0}^{T} \beta_{\frac {k'}{n}} (s')\beta_{\frac {k}{n}}(s)  \gamma_{t_{k'}, s'} (u')\gamma_{t_k, s} (u) dudu' \mu(dsds').
\end{eqnarray*}
Now by a  change of variables from $(u,u', s,s')$ to $ \frac{T}{n}(u,u', s,s') $, we obtain
\begin{eqnarray*}
 \mE
[{ {z}}_n(t)^2 ]
 &=& \left(\frac{T}{n}\right)^{2H+2}  \sum_{k,k' =0}^{ \lfloor \frac{nt}{T} \rfloor } \int_{k'}^{k'+1} \int_{k}^{k+1} \int_{0}^{n} \int_{0}^{n}   \varphi_{ {k'}, s'} (u')\varphi_{ k, s} (u) dudu' \mu(dsds')
 \\
 &=&   \left(\frac{T}{n}\right)^{2H+2}  \sum_{k,k' =0}^{ \lfloor \frac{nt}{T} \rfloor }\sum_{i,j=0,1} (-1)^{i+j}\tilde{e}_{ij},
\end{eqnarray*}
where 
\begin{eqnarray*}
\tilde{e}_{ij} &=& \int_{k'}^{k'+1} \int_{k}^{k+1} \int_{0}^{n} \int_{0}^{n}   \varphi_{ {k'}, s'}^{i} (u')\varphi_{ k, s}^{j} (u) dudu' \mu(dsds'),
\end{eqnarray*}
and $\varphi_{k,s}^{0} =\mathbf{1}_{[k,s]}$, $\varphi_{k,s}^{1} =\mathbf{1}_{[s, k+1]}$, $\varphi_{k,s} =\varphi_{k,s}^{0} -\varphi_{k,s}^{1} $, as defined as in the previous subsection. 
It is easy to see  that
\begin{eqnarray}\label{e3.25i}
\tilde{e}_{00} = \wt{Q}(k-k'), \quad\quad \text{and}\quad
\tilde{e}_{10} = \wt{R}(k-k').
\end{eqnarray}
  By a change of   variables from $(s,s')$ to $(k+1-s, k+1-s')$, we obtain 
  \begin{eqnarray}\label{e3.25}
\tilde{e}_{11}&=&
 \int_{ k-k'}^{ {k-k'+ 1} }
\int^{ {1  } }_{ {0} } (s'-(k-k')) s  \mu( ds ds'   )
=  \wt{Q}(k-k'),
  \end{eqnarray}
  where the second equation follows by exchanging the orders of the two integrals. 
 By changing the variables from $(s,s')$ to $(k'+1-s, k'+1-s')$  for $\tilde{e}_{11}$, we obtain 
 \begin{eqnarray*}
   \tilde{e}_{01} &=&
   \int_{0}^{1} \int_{k'-k}^{k'-k+1} (1-s') ( s-(k'-k)) \mu(dsds')
=  \wt{R}(k'-k),
\end{eqnarray*}
and so 
\begin{eqnarray}\label{e3.26}
  \sum_{k ,k'=0}^{ \lfloor \frac{nt}{T} \rfloor }  \tilde{e}_{01} =   \sum_{k ,k'=0}^{ \lfloor \frac{nt}{T} \rfloor }   \wt{R}(k'-k)
=   \sum_{k ,k'=0}^{ \lfloor \frac{nt}{T} \rfloor }   \wt{R}(k-k').
\end{eqnarray}
 In summary from \eref{e3.25i},  \eref{e3.25} and \eref{e3.26},  we obtain
 \begin{eqnarray}\label{eqn 2.4}
 \mE
[{ {z}}_n(t)^2 ]
   &= &
     2   \left(\frac{T}{n}\right)^{2H+2}
   \sum_{k,k' =0}^{ \lfloor \frac{nt}{T} \rfloor } 
  \left(
  \wt{Q  }(k-k') -\wt{R}(k-k') 
  \right)
  \nonumber
    \\
   & = &
  2  \left(\frac{T}{n}\right)^{2H+2}  
  \left(
  \sum_{p=0 }^{\lfloor \frac{nt}{T} \rfloor }
  \sum_{k' =  0}^{ \lfloor \frac{nt}{T} \rfloor  -p
    }
  ( \wt{Q}(p) - \wt{R}(p) )
    +
  \sum_{p=- \lfloor \frac{nt}{T} \rfloor  }^{ -1 }
   \sum_{k' =   -p }^{
    \lfloor \frac{nt}{T} \rfloor
    }
   ( \wt{Q}(p) - \wt{R}(p) )
   \right)
   \nonumber \\
  &:= &
    {\wt{q}}_1 +{\wt{q}}_2  .
   \end{eqnarray}

\noindent \textit{Step 2.}\quad In this step, we show the inequality \eref{eq 3.17}. 
Since $|\wt{Q} (p) - \wt{R}(p)|\sim p^{2H-3}$ for sufficiently large $p$, 
it is easy to see   that the series
$ \sum_{p \in \ZZ}  | \wt{Q}(p) - \wt{R}(p) | $ is convergent.  So we have the estimates
\begin{eqnarray}\label{e3.31}
{\tilde{q}}_1 &\leq& 2   \left(\frac{T}{n}\right)^{2H+2}   (\frac{nt}{T}+1) \sum_{p=0}^{\infty} |\wt{Q}(p) - \wt{R}(p)|
\end{eqnarray}
and
\begin{eqnarray}\label{e3.32}
{\tilde{q}}_2 &\leq&
 2   \left(\frac{T}{n}\right)^{2H+2}  (  \frac{nt}{T} +1) \sum_{p=-\infty}^{ -1 } |\wt{Q}(p) - \wt{R}(p)|.
\end{eqnarray}
 Applying  \eref{e3.31} and \eref{e3.32}  to \eref{eqn 2.4} yields
\begin{eqnarray}\label{e 3.5}
 n^{2H+1} \mE(  {z}_n(t)^2 ) &\leq& K( t+\frac{T}{n}).
 \end{eqnarray}
Take $s,t \in \Pi$. By replacing $t $ in \eref{e 3.5} by $t-s -\frac{T}{n} $ and noticing that $z_{n}(t)-z_{n}(s)$ and $z_{n}(t-s-\frac{T}{n})$ have the same distribution, we obtain 
\begin{eqnarray}\label{e 3.12}
n^{2H+1} \mE(|  {z}_n(t) -  {z}_n(s)|^2 )
&=& n^{2H+1} \mE(|  {z}_n(t -s-\frac{T}{n})|^2 )
\nonumber \\
&\leq& K(t-s).
\end{eqnarray}
  This completes the proof of \eref{eq 3.17}.

\medskip

 \noindent \textit{Step 3.}\quad 
 In this step, we show  the convergence of the process $(n^{H+\frac12} z_{n}, B)$. 
 Note that the finite dimensional distributions of $(n^{H+\frac{1}{2}}z_{n}, B)$ are Gaussian, so to show their convergences it suffices to show the convergences of their covariances.   We first consider the convergence of $n^{2H+1}\mE(|z_{n}(t)|^{2})$.  To   this aim, we write
    \begin{eqnarray}\label{e3.35}
{\tilde{q}}_1&=& 
     2    \left(\frac{T}{n}\right)^{2H+2}
  \sum_{p=0 }^{\lfloor \frac{nt}{T} \rfloor }
  (\lfloor \frac{nt}{T} \rfloor  -p+1)
  ( \wt{Q}(p) - \wt{R}(p) )
  \nonumber
    \\
    &=& 
     2  \left(\frac{T}{n}\right)^{2H+2}
\left(
\lfloor \frac{nt}{T} \rfloor
  \sum_{p=0 }^{\lfloor \frac{nt}{T} \rfloor }
          ( \wt{Q}(p) - \wt{R}(p) )
    -
     \sum_{p=0 }^{\lfloor \frac{nt}{T} \rfloor }
  (p-1)
    ( \wt{Q}(p) - \wt{R}(p) )
\right)
\nonumber
\\
&:=&  \wt{q}_{11} + \wt{q}_{12}.
   \end{eqnarray}
 First, it is easy to verify the following convergence:  
     \begin{eqnarray}\label{e3.36}
    \lim_{n \rightarrow \infty} n^{2H+1} {\wt{q}}_{11}
   &  =  &
     \lim_{n \rightarrow \infty} 2 n^{2H+1}     \left(\frac{T}{n}\right)^{2H+2}
 \lfloor \frac{nt}{T} \rfloor
  \sum_{p=0 }^{\lfloor \frac{nt}{T} \rfloor }
         ( \wt{Q}(p) - \wt{R}(p) )
        \nonumber
         \\
        &  =&
        2 T^{2H+1}   {t}    \sum_{p=0 }^{\infty }
         ( \wt{Q}(p) - \wt{R}(p) ) .
    \end{eqnarray}
On the other hand, since  
    $ | \sum_{p=0 }^{\lfloor \frac{nt}{T} \rfloor }
  (p-1)
    ( \wt{Q}(p) - \wt{R}(p) ) | \leq K n^{2H-1}  $, we have the convergence:
    \begin{eqnarray}\label{e3.37}
 \lim_{n \rightarrow \infty} n^{2H+1} {\wt{q}}_{12} &=& 0.
 \end{eqnarray}
Putting
together \eref{e3.36} and \eref{e3.37}, and taking into account \eref{e3.35}, we obtain:
 \begin{eqnarray}\label{eqn 3.7}
 \lim_{n \rightarrow \infty} n^{2H+1} {\wt{q}}_1 &=&   2   {T}^{2H+1}{t}   \sum_{p=0 }^{\infty }
         ( \wt{Q}(p) - \wt{R}(p) ) .
 \end{eqnarray}
 The quantity $\tilde{q}_{2}$ can be considered in a similar way. We can show that 
     \begin{eqnarray}\label{eqn 2.7}
    \lim_{n\rightarrow \infty} n^{2H+1} {\wt{q}}_2 &=& 2      {T}^{2H+1}{t}  \sum^{ -1}_{p= -\infty} ( \wt{Q}(p) - \wt{R}(p) )  .
    \end{eqnarray}
   Applying \eref{eqn 3.7} and \eref{eqn 2.7}
 to \eref{eqn 2.4}  we obtain
 \begin{eqnarray}\label{eqn 2.8}
  \lim_{ n \rightarrow \infty} n^{2H+1} \mE ({ {z}}_n(t)^2  ) =  2     {T}^{2H+1} \varrho {t}.
  \end{eqnarray}
  Take $s,t \in [0,T]$. By the same argument as in \eref{e 3.8i} and with the help of \eref{e 3.5} and \eref{eqn 2.8}, we can show that 
  \begin{eqnarray}\label{e3.34i}
  \lim_{ n \rightarrow \infty} n^{2H+1} \mE ( {z}_n(t)  {z}_n(s)  )
&=& 
     2T^{2H+1} \varrho    { ( t\wedge s)}{ }  .
  \end{eqnarray}
  
  On the other hand, by some elementary computation  (see Section \ref{section 6.2}), one can   show that 
\begin{eqnarray}\label{e3.34}
\lim_{n\rightarrow \infty } \mE(z_{n}(t) B_{r}) &=& 0.
\end{eqnarray}
Therefore, combining \eref{e3.34i} and \eref{e3.34}, 
we conclude        that the covariances of the finite dimensional distributions of    $ (n^{H+\frac12}  {z}_n , B )$ converge   to those of $(\sqrt{ {2 \varrho}{ } }T^{H+\frac12}   W , B )$. The proof is now complete.  \hfill$\Box$

\section{The strong convergence}\label{sec4}

We recall that  $X$  is  the solution of   equation \eref{e.1.1} and  $X^n$ is the continuous time interpolation of the Crank-Nicolson scheme defined in   \eref{4.1}.  In this section we   prove Theorem \ref{thm1.1} and some auxiliary results.

\medskip
\noindent \textit{Proof of Theorem \ref{thm1.1}:}\quad The proof is divided  into six  steps.
 
\noindent \textit{Step 1: Decomposition of the error process.}\quad In this step, we  derive a decomposition for the error process $Y_t : = X_t - X^{n}_t $\,, \,$t\in [0, T]$. 
For convenience  we set $\eta(t) =t_{k}$  for $t \in [t_{k}, t_{k+1})$ and $\ep(t) = t_{k+1}$  for $t \in (t_{k}, t_{k+1}]$.   Putting together equations \eref{e.1.1} and    \eref{4.1}, it is easily seen  that
  \begin{eqnarray}
Y_t  &= &    \int_0^t \left[ V  (X_s) -  V  (X^n_s)  \right] dB_s +
\frac12   \int_0^t \left[ V  (X^n_s)   - V  (X^n_{\eta(s)})\right] dB_s
\nonumber
\\
 &&
+
\frac12   \int_0^t \left[  V  (X^n_s) -    V  (X^n_{  \ep(s)   }) \right] dB_s  
\nonumber
\\
&=  & 
\sum_{j=0}^m  \sum_{i =1}^d \int_0^t  V_{ji}(s)  Y^{i}_s dB^j_s
+\frac12 J_1  (t)+\frac12 J_2   (t)  % \left[  V^i (X^n_{\ep(t)}) -    V^i (X^n_{   t  }) \right]  B_{\eta(t), t}
\,,
\label{eqn 3.1ii}
 \end{eqnarray}
where we have set for $t \in [0,T]$:
\begin{eqnarray*}
V_{ji}(s) &= &
 \int_0^1 \partial_{i} V_j (\theta X_s +(1-\theta) X^n_s ) d\theta,
\end{eqnarray*}
\begin{eqnarray*}
J_{1}(t) =  \int_0^t \left[ V  (X^n_s)   - V  (X^n_{\eta(s)})\right] dB_s\,,
 \quad\quad J_{2}(t) =   \int_0^t \left[  V  (X^n_s) -    V  (X^n_{  \ep(s)   }) \right] dB_s  \,,
\end{eqnarray*} 
and   we denote by    $\partial_{i}$   the partial differential operator with respect to the $i$th variable, that is, $\partial_{i}f(x)=\frac{\partial f}{\partial x_{i}}(x)$ for $f \in C^{1}$.
 In addition, 
   the chain rule for the Young integral enable us to write
\begin{eqnarray*}
 V  (X^n_s)   - V (X^n_{\eta(s)}) 
&=&\sum_{i=1}^{d}  \partial_{i} V  (X^{n}_{\eta(s)}) (X^{n,i}_{ s} - X^{n,i}_{\eta(s)} )   \\
&&+  \sum_{i,i'=1}^{d}
\int_{\eta(s)}^{s}  \int_{\eta(s)}^{u}   \partial_{i'} \partial_{i} V  (X^{n}_{v}) dX^{n,i'}_{v} dX^{n,i}_{u}\,.
\end{eqnarray*}
  Substituting the above expression into $J_1 (t)$, we obtain the following decomposition for $J_{1}(t)$
  \begin{eqnarray}\label{eq4.6}
    J_1  (t)
   &= &
   R_0   (t) + R_1  (t)
   , \quad t\in [0,T],
 \end{eqnarray}
where we define
\begin{eqnarray}\label{e4.5ii}
 R_1   (t) &= &    \int_0^t\left[
 \sum_{i,i'=1}^d \int^s_{ \eta(s)  } \int^{u}_{\eta(s) } \partial_{i'}\partial_{i} V  (X^n_{v}) dX^{n,i'}_{v} dX^{n,i}_{u} \right]
 dB_s\,
  \end{eqnarray}
and
\begin{eqnarray}
R_0   (t)
&= &    \int_0^t\left[
\sum_{i=1}^{d}  \partial_{i} V  (X^{n}_{\eta(s)})  (X^{n,i}_{ s} - X^{n,i}_{\eta(s)} ) 
 \right]
 dB_s
 \nonumber
\\
&= &
\frac12 \sum_{i=1}^{d}  \sum_{j,j'=0}^{m}  \int_0^t 
   \partial_{i}  V_{j}  (X^n_{ \eta(s)   })
    \left[ V_{j'}^{i}  (X^n_{   \ep(s) }) + V_{j'}^{i}  (X^n_{\eta(s )}) \right]    \int_{\eta(s)}^{s} dB^{j'}_{u}
     dB^{j}_s
  \,,
  \label{eqn3.1 i}
 \end{eqnarray}
 and  in the second equation of \eref{eqn3.1 i} we have used relation \eref{4.1}.

We can proceed similarly as in \eref{eq4.6} to derive  the corresponding decomposition for $J_{2}(t)$
 \begin{eqnarray*}
  J_2  (t)
 &= & - \wt{R}_0   (t) + \wt{R}_1   (t),\quad t\in [0,T],
 \end{eqnarray*}
where  
 \begin{eqnarray}\label{e4.5i}
 \wt{R}_1 (t) &= &      \int_0^t\left[
 \sum_{i,i'=1}^d \int_s^{  \ep(s)   } \int_{u}^{  \ep(s)    } \partial_{i'}\partial_{i} V  (X^n_{v}) dX^{n,i'}_{v} dX^{n,i}_{u} \right]
 dB_s\,,
 \\
\wt{R}_0  (t)
   &= &
\frac12 \sum_{i=1}^{d}  \sum_{j,j'=0}^{m}  \int_0^t 
   \partial_{i}  V_{j}  (X^n_{ \ep(s)   })
    \left[ V_{j'}^{i}  (X^n_{   \ep(s) }) + V_{j'}^{i}  (X^n_{\eta(s )}) \right]    \int_{s}^{\ep(s)} dB^{j'}_{u}
     dB^{j}_s  \,.
    \nonumber
   \end{eqnarray}

To further decompose the process $J_{1}$ and $J_{2}$, we introduce   the processes  $I_{1} $ and $I_{2}  $ defined on  $ \Pi$. Namely, for  $t \in \Pi\setminus \{ 0 \}$ we define
  \begin{eqnarray}
  I_{ 1}   (t) &= &  \sum_{ j,j'=0}^m
   \sum_{k=0}^{    {nt}/{T}   -1}
 (  \partial  V_{j}  V_{j'}) (X^n_{t_k})
   \int_{t_k}^{t_{k+1} }
\int_{t_{k}}^{s} d  B^{j'}_{u} 
     dB^{j}_s\,,
     \label{eqn 3.3 i}
     \\
     I_{ 2}  (t) &= & \sum_{ j,j'=0}^m
 \sum_{k=0}^{    {nt}/{T}   -1}
  ( \partial  V_{j} V_{j'} ) (X^n_{t_k})
   \int_{t_k}^{t_{k+1}   }  \int_{s}^{t_{k+1}} d B^{j'}_{u} dB^{j}_s
     \,,
     \label{eqn 3.3 ii}
     \end{eqnarray}
    and for  $t=0$ we set $I_{1}(0) = I_{2}(0) = 0$, where we  used the notation  $\partial =(\partial_{1},\dots, \partial_{d})$ and
    $\partial V_jV_{j'} = \sum_{i=1}^d \partial _i V_jV_{j'}$.
     Subtracting \eref{eqn 3.3 ii} from  \eref{eqn 3.3 i}  we obtain the following ``L\'evy area term''
     \begin{eqnarray}\label{eq4.12}
I_{1}(t) - I_{2}(t) =E_{1}(t) :=  \sum_{ j,j'=0}^m
 \sum_{k=0}^{    {nt}/{T}   -1}
  ( \partial  V_{j} V_{j'} ) (X^n_{t_k}) \zeta^{j'j}_{t_{k},t_{k+1}},
\end{eqnarray}
where we have denoted
\begin{eqnarray*}
\zeta_{s  t}^{ij} &=&
 \int_{s}^{  t }
 \int_{s}^{u}  dB^{i}_{v}   
     dB^j_u -  \int_{s}^{   t  } \int_{u}^{t} dB^{i}_{v}  dB^j_u, \quad
0\leq s\leq t\leq T.
\end{eqnarray*}  
Note that
a simple application of  Fubini's theorem to $ \zeta_{s, t}^{ij}  $ yields the identity
$\zeta_{s, t}^{ij} = - \zeta_{s, t}^{ji} $\,. So     expression \eref{eq4.12} can be    reduced to
 \begin{eqnarray}
E_{1}  (t)     &= & \sum_{j' < j }
   \sum_{k=0}^{   {nt}/{T} -1 } 
   \phi_{jj'} (X^n_{t_k})
    \zeta_{t_k, t_{k+1}   }^{j'j}
    \,,  \quad\quad t \in \Pi \, ,
     \label{eqn 3.4i}
      \end{eqnarray}
where $\phi_{jj'}$ is defined as $\phi_{jj'}=\partial  V_j V_{j'}
   -\partial  V_{j'} V_{j }$.
    In particular,     when    the driving process $B$ has dimension   one we are left with
      $E_{1}\equiv 0$.

With these calculations in hand, we can now decompose $J_{1}(t)+J_{2}(t)$ for  $t \in \Pi$ as follows: 
\begin{eqnarray}
J_{1}(t) + J_{2}(t)  & = &
      \left( I_1 (t) - I_2 (t) \right)
     +
      \left( R_0 (t) - I_1 (t) \right)
      +\left( I_2 (t)-\wt{R}_{0} (t) \right)  +R_1(t) +\wt{R}_1(t) 
    \nonumber
      \\
      &:= &
      E_1 (t) +E_2 (t) +E_3 (t) +E_4 (t) +E_5 (t)   \,.
      \label{eq4.5}
\end{eqnarray}

%$$
%\frac12 \sum_{i' =1}^d  \sum_{j,j'=0}^m %\sum_{k=0}^{\lfloor \frac{nt}{T} \rfloor} %\partial_{i'} V^i_j (X^n_{ t_k   })
%    [ V^{i'}_{j'} (X^n_{t_{k+1} \wedge t}) + %V^{i'}_{j'} (X^n_{t_k})]
%    \int_{t_k}^{t_{k+1}\wedge t}
%  \int^s_{t_k}  dB^{j'}_{s'}
%  dB^j_s
%$$

\

\noindent \textit{Step 2: Upper-bound for the Crank-Nicolson scheme.}\quad
It follows from Lemma 8.4 in  \cite{HLN2} that
  there exists a   constant $K$ such that
\begin{eqnarray}\label{eqn 3.0}
   \|X^n\|_{\infty}\vee \|X^n\|_{\beta}  & \leq & K +{K \|B\|_{\beta}^{1/\beta} }  .
\end{eqnarray}
Furthermore,   there exist  constants $K_0 $ and $K'_0$ independent of $n$ such that  for $0 \le s<t\le T$ and  $(t-s)^\beta  \|B\|_{\beta} \le K_0 $, we have
 \begin{eqnarray}\label{eqn 4.2}
 \| X^n \|_{s, t , \beta} &\leq &     K'_0  \|B \|_{\beta}  \,  .
 \end{eqnarray}

\noindent \textit{Step 3: Estimates of $E_{e}, 1\le e\le 5$. }\quad Take $ s,t \in \Pi$ such that $ s\leq t$. In this   step, we derive  a $L^{p}$-estimate of $ E_{e} (t) - E_{e}(s)$ for $e=1,\dots, 5$. We first show that for $e=2,3,4,5$ we have
\begin{eqnarray}\label{eq4.11}
\|E_{e}(t) - E_{e}(s)  \|_{p}& \leq & K n^{-2H}(t-s)^{\frac12}, \quad s,t \in \Pi,
\end{eqnarray}
where recall that $\| \cdot \|_{p}$ denotes the $L^p$-norm. 

 Let us start by bounding the term $E_2(t) -E_{2}(s) $,  $s,t \in \Pi$. Subtracting \eref{eqn 3.3 i}  from  \eref{eqn3.1 i} we obtain  
\begin{eqnarray}
E_{2}(t) &= &
 \sum_{i=1}^{d} \sum_{j,j'=0}^{m} \sum_{k=0}^{\frac{nt}{T}-1} \frac{1}{2} \partial_{i}V_{j} (X^{n}_{t_{k}}) \left[ V^{i}_{j'}(X^{n}_{t_{k+1}} )- V^{i}_{j'}(X^{n}_{t_{k }} ) \right]  \int_{t_{k}}^{t_{k+1}} \int_{t_{k}}^{s} dB^{j'}_{u} dB^{j}_{s}
\nonumber
 \\
 &=&
  \frac14 \sum_{i=1}^{d}   \sum_{j,j',j''=0}^m \sum_{k=0}^{  \frac{nt}{T} -1}  
\partial_{i}  V_j (X^n_{ t_k   })
    \int_{t_{k}}^{t_{k+1}} \partial V^{i}_{j'} (X^n_{v })   \left[ V_{j''}(X^{n}_{t_{k+1} } ) + V_{j''}(X^{n}_{t_{k}})  \right]  d B^{j''}_{v}
\nonumber    \\
     && 
   \quad\quad\quad\quad  \quad\quad\quad  \times
\int_{t_k}^{t_{k+1} }  \int^s_{t_k}  dB^{j'}_{u}
  dB^j_s  
  \,,
  \label{eq4.10}
   \end{eqnarray}
   where the second equation follows by applying the chain rule to $V^{i}_{j'}(X^{n}_{t_{k+1}} )- V^{i}_{j'}(X^{n}_{t_{k }} )$ and taking into account   equation \eref{4.1} for $X^{n}$. 
   Take $h^{n}_{t_{k}} = \partial_{i} V_{j} (X^{n}_{t_{k}})  \left[ V_{j''}(X^{n}_{t_{k+1} } ) + V_{j''}(X^{n}_{t_{k}})  \right] $ and $f_{v} = \partial V^{i}_{j'} (X^n_{v })  $. The  above expression becomes
   \begin{eqnarray}\label{e4.15}
E_{2}(t) &= & \frac{1}{4}  \sum_{i=1}^{d} \sum_{j,j',j''=0}^{m} \sum_{k=0}^{\frac{nt}{T}-1}    \int_{t_{k}}^{t_{k+1}} 
\int_{t_k}^{t_{k+1} }  \int^s_{t_k}    f_{v} h_{t_{k}}^{n}  dB^{j'}_{u}
  dB^j_s    d B^{j''}_{v}.
\end{eqnarray}
It is easy to verify that  the triple integral on the right-hand side of \eref{e4.15} is equal to 
 \begin{eqnarray*}
 && 
  \int_{t_{k}}^{t_{k+1}} 
\int_{t_k}^{v }  \int^s_{t_k}    f_{v} h_{t_{k}}^{n}  dB^{j'}_{u}
  dB^j_s    d B^{j''}_{v} 
 +  
    \int_{t_{k}}^{t_{k+1}} 
\int_{t_{k}}^{s }  \int^v_{t_k}    f_{v} h_{t_{k}}^{n}  dB^{j'}_{u}  d B^{j''}_{v}
  dB^j_s 
  \\
  && 
  + \int_{t_{k}}^{t_{k+1}} 
\int_{t_{k}}^{s }  \int^{u}_{t_{k}}    f_{v} h_{t_{k}}^{n}   d B^{j''}_{v}  dB^{j'}_{u}
  dB^j_s   
\,   .
\end{eqnarray*}
Substituting the above expression  into \eref{e4.15}, we obtain an expression of  $E_{2}(t)$ of the form of \eref{e6.23}.  One can show, with the help of the estimate  of $X^{n}$ in \eref{eqn 3.0}, that $f$ and $h^{n}$ satisfy the conditions in  Lemma \ref{lem5.2}.  So applying Lemma \ref{lem5.2} to  $E_{2}(t)$  we obtain the estimate
 \eref{eq4.11}
for $e=2$. 
Estimate \eref{eq4.11} still holds true for the cases when $e=3,4,5$. The proof is   based on Lemma \ref{lem5.2} and is similar to the case $e=2$. We omit the details. This completes the proof of \eref{eq4.11}.
 
 Now we consider the process  $E_1(t)$,  $t \in \Pi  $. To this aim, we consider the decomposition
 \begin{eqnarray}\label{e4.21}
E_{1} (t) &= &   \sum_{0 \neq j' < j }
   \sum_{k=0}^{   {nt}/{T} -1 } 
   \phi_{jj'} (X^n_{t_k})
    \zeta_{t_k, t_{k+1}  }^{j'j} + \sum_{0= j' < j }
   \sum_{k=0}^{   {nt}/{T}-1  } 
   \phi_{jj'} (X^n_{t_k})
    \zeta_{t_k, t_{k+1}   }^{j'j}
    \nonumber
    \\
   & := & E_{11}(t)+  E_{12} (t).
\end{eqnarray}
Expression \eref{e4.21} and Lemma  \ref{lem 2.4} together suggest   to consider the following ``weight-free'' random sum corresponding to $E_{11}$	 
\begin{eqnarray*}
{g}_{n} (t) &= &
n^{2H-\frac 12 } \sum_{0 \neq j' < j }
   \sum_{k=0}^{\lfloor \frac{nt}{T} \rfloor} 
     \zeta_{t_k, t_{k+1}   }^{j'j} .
\end{eqnarray*}
     It follows from relation \eref{e 3.1i}  in    Proposition \ref{prop 3.1}  that   $ {g}_{n}$ satisfies the assumptions in   Lemma \ref{lem 2.4}.
     Indeed, by Proposition \ref{prop 3.1} the following estimate holds true for all $s,t\in \Pi$
     \begin{eqnarray}\label{eqn4.22}
 \mE( [  {g}_n(t)-  {g}_n(s) ]^2 )^{\frac12} &\leq& K |t-s|^{\frac12}.
\end{eqnarray}
Furthermore, since $g_{n}(t)-g_{n}(s)$ is a random variable in the second chaos of $B$, by an hyper-contractivity argument we can  show that   estimate \eref{eqn4.22}   holds in the $L^p$-norm for all $p\geq 1$.
       Take $f =   \phi_{jj'} (X^{n}_{\cdot})$, $\beta'=\frac12$, $ \frac12< \beta<H$, $p=p'=q'=2$. Then applying  Lemma \ref{lem 2.4} to $E_{11} $   we obtain the estimate
 \begin{eqnarray}\label{eqn 3.7i}
\| E_{11}(t)-E_{11} (s) \|_p  & \leq   & Kn^{-2H+\frac12} (t-s)^{\frac12}, \quad s,t \in \Pi.
\end{eqnarray}

We proceed similarly to show the estimate for $E_{12}$.   We first define the ``weight-free'' random sum corresponding to $E_{12} (t)$
\begin{eqnarray*}
  \tilde{g}_{n} (t) &= & n^{1/2+H} \sum_{0 = j' < j }
   \sum_{k=0}^{\lfloor \frac{nt}{T} \rfloor} 
     \zeta_{t_k, t_{k+1}   }^{j'j}.
\end{eqnarray*}
Then  as in \eref{eqn4.22},   estimate     \eref{eq 3.17} in Proposition \ref{prop 3.2} together with  some hyper-contractivity arguments yields that  
$ \tilde{g}_{n} $ 
  satisfies the conditions in Lemma \ref{lem 2.4} for $\beta'=\frac12$ and $p=q'=2$. Taking $\frac12<\beta<H$, $q'=2$ and $f =   \phi_{jj'} (X^{n}_{\cdot})$ as before  and applying Lemma \ref{lem 2.4} to $E_{12}$,  we obtain the estimate
 \begin{eqnarray}\label{eq 4.13} 
\| E_{12}(t)-E_{12} (s) \|_p & \leq & Kn^{-H-\frac12} (t-s)^{\frac12}, \quad s,t \in \Pi.
\end{eqnarray}

In summary of     relations  \eref{eq4.11},  \eref{eqn 3.7i} and \eref{eq 4.13},  and taking into account the fact that $E_{11}=0$ when $m=1$ and $E_{11}=E_{12}=0$ when $m=1$ and $V_{0} \equiv 0$, we obtain
\begin{eqnarray}\label{eq4.13}
\sum_{e=1}^{5} \|E_{e}(t) - E_{e}(s) \|_{p} & \leq & 
  K  (t-s)^{\frac12} /\vartheta_{n} ,  \quad s,t \in \Pi.
  \end{eqnarray}

 \noindent \textit{Step 4: Upper-bounds for the Jacobian.}\quad In this step, we consider some linear equations associated with $X^{n}$ and $X$. 
  Let $\Lambda^{n }  =  \left( \Lambda^{n,i}_{i'}  \right)_{1\leq i,i'\leq d}$ be the solution of the linear equation
\begin{eqnarray}
\Lambda^{n,i}_{{i'}   }(t) &= & \delta^i_{i'}+ \sum_{j=0}^m \sum_{i'' =1}^d\int_0^t V^i_{j  i''} (s) \Lambda^{n,i''}_{{i'}   } (s) dB^{j }_s, ~~~~~~~i,i'=1, \dots, d,\quad t \in [0, T].
\label{eqn 3.1i}
\end{eqnarray}
Here $\delta^i_{i'} $ is the Kronecker function, that is, $\delta^i_{i'} =1$ when $i= i'$ and $\delta^i_{i'} = 0 $ otherwise. The $d\times d$ matrix $\Lambda^n (t) $ is invertible. We denote  its inverse  by $\Gamma^n(t)   $. It is easy to verify that
$\Gamma^n$ satisfies the equation
\begin{eqnarray*}
\Gamma^{n,i}_{{i'}   }(t) &=& \delta^i_{i'} - \sum_{j=0 }^m \sum_{i'' =1}^d\int_0^t \Gamma^{n,i }_{{i'' }   } (s) V^{i''}_{j  i'} (s)  dB^{j }_s, ~~~~~~~i,i'=1, \dots, d,\quad t \in [0 , T].
\end{eqnarray*}
By the product rule of Young integrals, and taking into account (\ref{eqn 3.1ii}), it is easy to verify that  
\begin{eqnarray}\label{eqn 3.5}
Y_t &= &
 \frac12  \Lambda^{n }_{  t}   \sum_{i=1}^{2}  \int_{0}^{ t }
    \Gamma^{n }_{s}
      dJ_{i}(s)
 , \quad t\in [0, T]\,.
      \end{eqnarray}
Applying Lemma 3.2 (ii) in \cite{HLN} and taking into account   \eref{eqn 4.2}, we  obtain the estimate
 \begin{eqnarray*} 
 \| \Lambda^n\|_{\infty} \vee \| \Lambda^n\|_{\beta} \vee \|  \Gamma^n \|_{\infty} \vee \|  \Gamma^n \|_{\beta} &\leq& K e^{K \|B\|_{\beta}^{1/\beta} }  .
\end{eqnarray*}
It follows from  Fernique's lemma that
for $p\geq 1$ we have
 \begin{eqnarray}\label{e4.22}
\left\| \| \Lambda^n\|_{\infty}\right\|_{p} \vee \| \| \Lambda^n\|_{\beta} \|_{p}\vee \|\|  \Gamma^n \|_{\infty}\|_{p} \vee \|\|  \Gamma^n \|_{\beta}\|_{p} &\leq& K   .
\end{eqnarray}

Let $\Lambda  = \left( \Lambda^i_{i'}  \right)_{1\leq i,i'\leq d}$  be the solution of the following equation,
\begin{eqnarray}\label{e 3.3}
\Lambda^{ i}_{{i'}   }(t) &= &
 \delta^i_{i'}+ \sum_{j=0}^m \sum_{i'' =1}^d\int_0^t
 \partial_{i''} V^i_j (  X_s  )\Lambda^{ i''}_{{i'}   } (s) dB^{j }_s, 
 \end{eqnarray}
for  $t \in [0, T]
$, $ i,i'=1, \dots, d$, and denote by $\Gamma(t)$ the inverse of $\Lambda(t)$. As before, we can show  that $\Gamma$ satisfies the equation
\begin{eqnarray*}
\Gamma^{ i}_{{i'}   }(t) &=& \delta^i_{i'} - \sum_{j=0 }^m \sum_{i'' =1}^d\int_0^t \Gamma^{n,i }_{{i'' }   } (s) \partial_{i'}V_{j}^{i''}(X_{s})   dB^{j }_s
\end{eqnarray*}
for   $t \in [0 , T]$, $ i,i'=1, \dots, d$.
 It follows from Lemma 3.1 in \cite{HLN} that the estimate \eref{e4.22} still holds true if we replace $\Lambda^{n}$ and $\Gamma^{n} $ in \eref{e4.22} by  $\Lambda$ and $\Gamma $. 

\medskip

\noindent \textit{Step 5: Estimates of \,$ \Gamma^{n} \, Y  $.} \quad     
In this step, we consider the   process $\Gamma^{n}  Y $.  
Multiplying both sides of \eref{eqn 3.5} by $\Gamma^{n}_{t}$, we obtain the   expression
 \begin{eqnarray*} 
\Gamma^{n}_{t}\, Y_{t} &=&  \frac12 \sum_{i=1}^{2}\int_{0}^{t } \Gamma^{n}_{u} d J_{i}(u).
\end{eqnarray*}
  By writing $ \Gamma^{n}_{u} = \Gamma^{n}_{\eta(u)} + (  \Gamma^{n}_{u} -  \Gamma^{n}_{\eta(u)}  ) $   we obtain the following decomposition for     $s,t \in \Pi$, $s\le t$
 \begin{eqnarray}\label{eq4.15}
  \sum_{i=1}^{2}\int_{s}^{t } \Gamma^{n}_{u} d J_{i}(u) &=  &  \sum_{i=1}^{2}\int_{s}^{t } \Gamma^{n}_{\eta(u)} d J_{i}(u) +  \sum_{i=1}^{2} \int_{s}^{t } \int_{\eta(u)}^{u} d \Gamma^{n}_{v} d J_{i}(u).
\end{eqnarray}
Revoking the decomposition \eref{eq4.5} we    get
 \begin{eqnarray}\label{e4.29i}
 \sum_{i=1}^{2}\int_{s}^{t } \Gamma^{n}_{\eta(u)} d J_{i}(u) 
 &=&
  \sum_{e=1}^{5} \sum_{t_{k}=s}^{t-\frac Tn}\Gamma^{n}_{t_{k}} (E_{e}(t_{k+1}) - E_{e} (t_{k})).
\end{eqnarray}
 For   simplicity, we will denote the right-hand side of \eref{e4.29i} as
\begin{eqnarray}\label{eq4.29}
   \sum_{t_{k}=s}^{t-\frac Tn}\Gamma^{n}_{t_{k}} (E_{e}(t_{k+1}) - E_{e} (t_{k})) 
 & :=& \int_{s}^{t } \Gamma^{n}_{\eta(u)} d E_{e}(u).
\end{eqnarray}
Note, however, that  equation \eref{eq4.29} is only valid   for $s,t \in \Pi$ since $E_{e}$, $e=1,\dots, 5$ are only  defined on $\Pi$. Now substituting \eref{e4.29i} into \eref{eq4.15} and taking into account \eref{eq4.29} we get
\begin{eqnarray*}%\label{e4.33}
 \sum_{i=1}^{2}\int_{s}^{t } \Gamma^{n}_{u} d J_{i}(u) &=  & \sum_{e=1}^{5}  \int_{s}^{t } \Gamma^{n}_{\eta(u)} d E_{e}(u) +  \sum_{i=1}^{2} \int_{s}^{t } \int_{\eta(u)}^{u} d \Gamma^{n}_{v} d J_{i}(u).
\end{eqnarray*}
As in \eref{eqn 3.7i}, we   handle the term $\int_{s}^{t } \Gamma^{n}_{\eta(u)} d E_{e}(u) $ on the right-hand side of \eref{e4.33}   by 
  Lemma \ref{lem 2.4}.   
    Take   $ \hat{g}_{n}(t) = \vartheta_{n} E_{e}(t) $, $t \in \Pi$ and $f  = \Gamma^{n} $, and let $\beta, \beta', p, p', q'$ be as before. Then 
   estimate \eref{eq4.13} shows that $\hat{g}_{n}$ satisfies the conditions in Lemma \ref{lem 2.4}. 
   Applying Lemma \ref{lem 2.4} to $\int_{s}^{t } \Gamma^{n}_{\eta(u)} d E_{e}(u) $ and invoking expression \eref{eq4.29} we obtain
\begin{eqnarray}\label{eq 4.17}
\Big\| \sum_{e=1}^{5} \int_{s}^{t } \Gamma^{n}_{\eta(u)} d E_{e}(u) \Big\|_{p} & \leq   & 
K (t-s)^{\frac12} / \vartheta_{n}\,, \quad  s,t \in \Pi.
\end{eqnarray}

We turn to the second term in \eref{eq4.15}. By the definition of $\Gamma^{n}$ and $J_{1}$ we have
\begin{eqnarray*}
 	\int_{t_{k}}^{t_{k+1} } \int_{t_{k}}^{u} d \Gamma^{n}_{v} d J_{1}(u)     & =  & 
  \frac12  \sum_{jj'=0 }^m \sum_{i,i' =1}^d	\int_{t_{k}}^{t_{k+1} } \int_{t_{k}}^{u} 
	\left( - \Gamma^{n }_{{i' }   } (v) V^{i'}_{j , i} (v) \right)  dB^{j }_v 
	\int_{t_{k}}^{u} \partial V_{j'}^{i} (X^n_{r}) d X^{n}_{r}    dB^{j'}_u  .
\end{eqnarray*}
One can show  that
\begin{eqnarray*}
  \int_{s}^{t } \int_{\eta(u)}^{u} d \Gamma^{n}_{v} d J_{1}(u)  
&=&
  \sum_{k=\frac{ns}{T}}^{\frac{nt}{T}-1} \int_{t_{k}}^{t_{k+1}} \int_{t_{k}}^{u} d \Gamma^{n}_{v} d J_{1}(u)   
\end{eqnarray*}
has  the form of \eref{e6.23}. 
Applying Lemma \ref{lem5.2} we obtain\begin{eqnarray}\label{eq 4.18}
\left\| \int_{s}^{t } \int_{\eta(u)}^{u} d \Gamma^{n}_{v} d J_{i}(u)  \right\|_{p}
 &\leq & K n^{-2H}(t-s)^{1/2}
\end{eqnarray}
for $i=1$. This estimate still holds true in the case  $i=2$, and the proof is  similar.
Substituting \eref{eq 4.17} and \eref{eq 4.18}   into \eref{eq4.15}  we obtain the estimate  
\begin{eqnarray}\label{eq 4.20}
\left\| \sum_{i=1}^{2 } \int_{s}^{t } \Gamma^{n}_{u} d J_{i}(u) \right\|_{p}  & \leq & 
K (t-s)^{\frac12} /\vartheta_{n} 
\end{eqnarray}
for $s,t \in \Pi$.

It is easy to see that
  \begin{eqnarray*} 
 \left\| \int_{t_{k}}^{t } \Gamma^{n}_{u}    d J_{e}(u) \right\|_{p} \leq\, & Kn^{-2H},  \quad t\in [t_{k}, t_{k+1}].
\end{eqnarray*}
Combining this estimate with \eref{eq 4.20} we obtain the inequality
 \begin{eqnarray}\label{eq 4.21}
\sup_{t\in [0,T]}
\left\|  \sum_{i=1}^{2}\int_{0}^{t } \Gamma^{n}_{u} d J_{i}(u)  \right\|_{p}  & \leq & 
 K / \vartheta_{n}   .
\end{eqnarray}

\medskip
  \noindent \textit{Step 6: Conclusion.}\quad 
  The   inequality \eref{eq4.2}  follows by applying the H\"older inequality to \eref{eqn 3.5} and using  the estimate \eref{eq 4.21} and  the estimate \eref{e4.22} for $ \Lambda^{n} $.     \hfill $\Box$
  
  \medskip
 In the last part of this section we will show some technical estimates that will be used in the proof of the convergence in law of the error. 
 
 \begin{lemma}
 Under the assumptions and notation of Theorem \ref{thm1.1},  the error process $Y=X-X^n$ satisfies the following relation for all $s,t\in \Pi$
 \begin{eqnarray}\label{eq4.3}
\mE \left(\left| Y_{t} - Y_{s} \right|^{p} \right)^{1/p} &\leq & 
  K    |t-s|^{\frac 12}/\vartheta_n.
   \end{eqnarray}
  \end{lemma}

  \noindent \textit{Proof: } \quad 
  Invoking the expression \eref{eqn 3.5} of $Y$,   we can write
 \begin{eqnarray}\label{e4.33}
Y_t-Y_{s} &= &
 \frac12 ( \Lambda^{n }_{  t}  -\Lambda^{n }_{  s}  ) \sum_{i=1}^{2}  \int_{0}^{ t }
    \Gamma^{n }_{u}
      dJ_{i}(u) + 
       \frac12  \Lambda^{n }_{  s}  \sum_{i=1}^{2}  \int_{s}^{ t }
    \Gamma^{n }_{u}
      dJ_{i}(u).
\end{eqnarray}
 The   inequality \eref{eq4.3} then follows by  applying the H\"older inequality to \eref{e4.33}
and by taking into account   the estimates     \eref{e4.22} and  \eref{eq 4.20}.   This completes the proof. 
\hfill $\Box$

\medskip

The following  lemma is a convergence result for the processes $\Lambda^{n}$ and $\Gamma^{n}$. 

 \begin{lemma}\label{lem4.4}
 Take $\be: \frac12 < \beta<H$.
 Let $\Lambda^{n}$ and $\Lambda$ be the solutions of equations \eref{eqn 3.1i} and \eref{e 3.3}, respectively, and let $\Gamma^{n}$ and $\Gamma$ be their inverses. Then we have 
 \begin{eqnarray}\label{eq4.24}
\| \Lambda^{n} -\Lambda \|_{\beta, p } +
  \| \Gamma^{n} -\Gamma \|_{\beta, p } 
 &\leq &  Kn^{1-2\beta}.
\end{eqnarray}
	\end{lemma}
\noindent \textit{Proof: } \quad See Section \ref{section 6.3i}.
\hfill $\Box$

\bigskip

We end this  section with      the following auxiliary results.
The reason we put these results here is because they are concerned with $\Gamma$.   As in \eref{eq4.29}, for the sake of conciseness we will denote   
\begin{eqnarray*}
  \sum_{t_{k}=s}^{t-\frac Tn}\Gamma^{n}_{t_{k}} (E_{11}(t_{k+1}) - E_{1} (t_{k})) 
 & :=& \int_{s}^{t } \Gamma^{n}_{\eta(u)} d E_{11}(u) \quad   \text{ for } s,t \in \Pi.
\end{eqnarray*}
The integral  $\int_{s}^{t} \Gamma^{n}_{\eta(u)} dE_{12} (u)$  is defined similarly. 
\begin{lemma}\label{lem4.2}
We continue to use 
  the notation of  in Theorem \ref{4.1}.  Let $s,t \in \Pi$, $s\le t$.  If  $m>1$, we have the estimate
\begin{eqnarray}\label{eq4.21}
 \sup_{s,t\in \Pi} \left\| \sum_{i=1}^{2 } \int_{s}^{t } \Gamma^{n}_{u} d J_{i}(u) - \int_{s}^{t} \Gamma^{n}_{\eta(u)} d E_{11}(u) \right\|_{p} & \leq  & 
K n^{-\frac12-H}.
\end{eqnarray}
In the case when $m=1$,   we have the estimate
\begin{eqnarray}\label{eq4.22}
  \sup_{s,t\in \Pi} \left\| \sum_{i=1}^{2 } \int_{s}^{t } \Gamma^{n}_{u} d J_{i}(u) - \int_{s}^{t} \Gamma^{n}_{\eta(u)} d E_{12}(u) \right\|_{p}  &\leq  & 
K n^{-2H}.
\end{eqnarray}
Take $\be: \frac12 < \beta<H$. Assume that $m=1$ and $V_{0}\equiv 0$.   Then  we have the estimate
\begin{eqnarray}\label{e4.23}
 \sup_{ t\in [0,T]} \left\| \sum_{i=1}^{2 } \int_{0}^{t } \Gamma^{n}_{u} d J_{i}(u) - \sum_{e=2}^{5 }\int_{0}^{\eta(t)} \Gamma^{n}_{\eta(u)} d E_{e}(u) \right\|_{p}  & \leq  & K_{\beta}n^{1-4\beta},
\end{eqnarray}
where $K_{\be} $ is a constant depending on $\be$.
\end{lemma}

 \noindent \textit{Proof: } \quad
By subtracting $\int_{s}^{t} \Gamma^{n}_{\eta(u)} d E_{11}(u) $ from both sides of \eref{eq4.15}
we obtain
\begin{eqnarray}\label{e4.27}
 &&\sum_{i=1}^{2}\int_{s}^{t } \Gamma^{n}_{u} d J_{i}(u) -\int_{s}^{t} \Gamma^{n}_{\eta(u)} d E_{11}(u)  
 \nonumber
 \\
 &&=\,  \int_{s}^{t} \Gamma^{n}_{\eta(u)} d E_{12}(u)  +  \sum_{e=2}^{5}\int_{s}^{t } \Gamma^{n}_{\eta(u)} d E_{e}(u) +  \sum_{i=1}^{2} \int_{s}^{t } \int_{\eta(u)}^{u} d \Gamma^{n}_{v} d J_{i}(u).
\end{eqnarray}
Similar to the proof of the estimate    \eref{eq 4.17}, we can show that the first   and second  terms on the right-hand side of \eref{e4.27} are bounded by $Kn^{\frac{1}{2}+H}$ and $K n^{-2H}$, respectively. On the other hand, we have shown    in \eref{eq 4.18} that the third term is bounded by $Kn^{-2H}$. In summary, we obtain      the estimate \eref{eq4.21}. 
The estimate \eref{eq4.22} can be shown in a similar way. 
 The proof of   estimate   \eref{e4.23} is included in   Section \ref{section 6.3}.  \hfill $\Box$

\section{Asymptotic error distribution}\label{sec5}
In this section, we prove Theorem \ref{thm1.2}.

\noindent \textit{Proof of Theorem \ref{thm1.2}:}\quad  The proof will be done in four  steps.

\medskip
 \noindent \textit{Step 1.}\quad  We first assume that $m>1$ or $V_{0}\not\equiv 0$. By Theorem 13.5 in  \cite{Bill} and taking into account  inequality \eref{eq4.3}, to prove the weak convergence of $( \vartheta_{n}( \tilde{X}  - \tilde{X}^{n}) , B)$  it suffices to show  the convergence of its   finite dimensional distributions (f.d.d.).  By \eref{eqn 3.5} we have
 \begin{eqnarray*} 
 \tilde{X}_{t} - \tilde{X}^{n}_{t} =
 X_{t_{k}} - X^{n}_{t_{k}}
 &= & \frac12 \Lambda^{n}_{ t_{k}}     \sum_{i=1}^{2}\int_{0}^{t_{k} } \Gamma^{n}_{u} d J_{i}(u) 
\end{eqnarray*}
  for $t\in [t_{k},t_{k+1})$. 
  
\medskip

 \noindent \textit{Step 2.}\quad
Assume that    $m>1$. Set
 \begin{eqnarray*}
S^{n}(t) &=&\frac12 \Lambda^{n}_{t_{k}}  \int_{0}^{t_{k}} \Gamma^{n}_{\eta(s)} d E_{11}  (s)
\end{eqnarray*}
for $t \in [t_{k},t_{k+1})$. 
 It follows from   the estimate \eref{eq4.21} in Lemma \ref{lem4.2} that the difference $\vartheta_{n}\|S^{n}(t) - ( \tilde{X}_{t} - \tilde{X}^{n}_{t} ) \|_{p}$ is uniformly bounded by $\vartheta_{n} n^{-\frac12 - H}$ and thus converges to zero as $n\rightarrow \infty$. This implies that  the limit of the finite dimensional distributions of   $  (\vartheta_{n}( \tilde{X}  - \tilde{X}^{n}) , B)$ is equal to that of   
$\left(  \vartheta_{n} S^{n}  , ~B  \right)$.

 Set
 \begin{eqnarray*}
S(t) &=& \frac12 \Lambda_{ t_{k}}  \int_{0}^{t_{k+1}} \Gamma _{\eta(s)} d E_{11}  (s)
%\\
%&=&  \Lambda_{ t_{l+1}}  \sum_{0 \neq j' < j } \sum_{k=0}^{ l }  \Gamma_{ t_{k}    } 
%    \phi_{jj'} (X_{t_k})
%     \zeta_{t_k, t_{k+1}   }^{j'j}
%     \\
%     &=&  \Lambda_{ t_{l+1}}  \sum_{0 \neq j' < j } \int_{0}^{t_{l+1}} \Gamma_{ \eta(s)    } 
%    \phi_{jj'} (X_{ \eta(s)  })
%     \zeta_{ \eta(s)  , s  }^{j'j}
\end{eqnarray*}
 for $t\in [t_{k} , t_{k+1})$.
Then    we have
 \begin{eqnarray}\label{e5.8}
S^{n}(t_{k}) -  S (t_{k}) 
%=\,&   \Lambda^{n}_{ t} \int_{0}^{t} \Gamma^{n}_{ \eta(s)} d E_{11}  (s)
%-  \Lambda_{ t}  \sum_{0 \neq j' < j } \sum_{k=0}^{  \frac{nt}{T} -1  }  \Gamma _{t_{k} }    
%    \phi_{jj'} (X_{t_k})
%     \zeta_{t_k, t_{k+1}   }^{j'j} 
%  \\
   &=&
   \frac12 \sum_{0 \neq j' < j }  \int_{0}^{t_{k}} \left[   \Lambda^{n}_{t_{k}} \Gamma^{n}_{ \eta(s) } \phi_{jj'} (X^{n}_{\eta(s)}) - \Lambda_{t_{k}} \Gamma_{ \eta(s)  } \phi_{jj'} (X_{\eta(s)})\right]  d \zeta^{j'j}_{\eta(s), s }
 \nonumber
   \\
   && - \frac12 \Lambda_{t_{k}}  \int_{t_{k}}^{t_{k+1}}  \Gamma_{\eta(s)} d E_{11}(s) . 
\end{eqnarray}
It is easy to see that the $L^p$-norm of the second term in the right-hand side of \eref{e5.8} is bounded by $Kn^{-2H}$. 
  On the other hand, with the help of    Lemma \ref{lem4.4}, one can show    that
\begin{eqnarray*}%\label{eq4.25}
\|  \Lambda^{n}_{t} \Gamma^{n}_{\cdot} \phi_{jj'} (X^{n}_{\cdot}) - \Lambda_{t} \Gamma_{\cdot}  \phi_{jj'} (X_{\cdot}) \|_{\beta, p } & \leq & K n^{1-2\beta}.
\end{eqnarray*}
So by taking $f = \Lambda^{n}_{t} \Gamma^{n}_{\cdot} \phi_{jj'} (X^{n}_{\cdot}) - \Lambda_{t} \Gamma_{\cdot}  \phi_{jj'} (X_{\cdot})$ and $\zeta_{k,n} = \zeta^{j'j}_{t_{k},t_{k+1}}$ in Lemma \ref{lem 2.4}, we see that   the first term in the right-hand side of \eref{e5.8} is bounded by $K n^{1-2\beta +1/2 -2H}$. In summary of these two estimates,  we obtain 
\begin{eqnarray*}
\left\| S^{n}(t) -  S(t) \right\|_{p} & \leq & K n^{1-2\beta +1/2 -2H} \vee n^{-2H}
\end{eqnarray*}
for $t \in \Pi$, and thus for $t \in [0,T]$. 
Therefore,   the f.d.d. convergence   of $\left(  \vartheta_{n} S^{n}  , B  \right)$ is the same as that of
$
\left(  \vartheta_{n}  S   , B \right)$. 

Applying Proposition \ref{prop2.8} to the process $
\left(  \vartheta_{n}  S   , B \right)$  and taking into account the weak convergence result in Proposition \ref{prop 3.1}, we obtain that the f.f.d. of $
\left(  \vartheta_{n}  S   , B \right)$ converge to that of 
$
\left(  U   , B \right)$, where 
\begin{eqnarray*} 
U_{t} &=&  T^{2H-\frac12}  \sqrt{  \frac{  \kappa}{ 2 } }   \Lambda_{t}  \sum_{1 \leq j' < j \leq m } \int_{0}^{t} \Gamma_{s} \phi_{jj'}(X_{s} ) dW^{j'j}_{s}   .
\end{eqnarray*}
  The convergence \eref{5.3} follows from the fact that $\{U_{t}, \, t\in [0,T]\}$  solves the SDE \eref{5.4}.

\medskip

 \noindent \textit{Step 3.}\quad
We assume $m=1$ and $V_{0} \not \equiv 0$. The  estimate \eref{eq4.22} implies that   the f.d.d. convergence   of $  (\vartheta_{n}(\tilde{X}  - \tilde{X}^{n}) , B)$ is equal to that of    $( \vartheta_{n}\wt{S}^{n}, B )$, where 
  \begin{eqnarray*}  
\wt{S}^{n}_{t} &=& \frac12  \Lambda^{n}_{ \eta(t)}   \int_{0}^{\eta(t) } \Gamma^{n}_{\eta(s)} d E_{12}  (s) .
\end{eqnarray*}
As in the case $m>1$, with the help of   Lemma \ref{lem4.4}     we can show   that 
  the convergence of the f.d.d. of  $( \vartheta_{n}\wt{S}^{n}, B )$ is the same as that of  $( \vartheta_{n}\wt{S} , B )$, where
\begin{eqnarray*}
\wt{S}_{t} &=& \frac12  \Lambda_{ \eta(t)}   \sum_{k=0}^{ \lfloor \frac{nt}{T} \rfloor  }  \Gamma_{ t_{k}    } 
    \phi_{10} (X_{t_k})
     \zeta_{t_k, t_{k+1}   }^{01}  .
\end{eqnarray*}
Applying Proposition \ref{prop2.8} to the   above process and taking into account the weak convergence result in Proposition \ref{prop 3.2}, we obtain that  its f.d.d. converges to those of  $(\wt{U}, B)$, where 
\begin{eqnarray*} 
\wt{U}_{t} &=&   T^{H+\frac12} \sqrt{\frac{\varrho}{2}} \Lambda_{t}  \int_{0}^{t} \Gamma_{s} \phi_{10} (X_s) d W_{s} 
\end{eqnarray*}
as $n\rightarrow \infty$. The convergence \eref{5.3} follows from the fact that $\{ \wt{U}_{t}, \, t\in [0,T]\}$ solves  equation \eref{5.5}.

\medskip

 \noindent \textit{Step 4.}\quad
We consider the case when  $m=1$ and $V_{0}\equiv 0$. The convergence \eref{5.6} is clear for $t=0$. In the following, we consider $t>0$.  
  The estimate \eref{e4.23} implies that the $L^p$-convergence of $n^{2H}(\tilde{X}_{t} - \tilde{X}^{n}_{t})$ is the same as   that of  
\begin{eqnarray}\label{e5.11ii}
  \frac12  n^{2H}  \Lambda_{  \eta(t)}^{n}   \sum_{e=2}^{5}  \int_{0}^{ \eta(t)  }
    \Gamma_{\eta(s)}^{n}
      dE_{e}(s).
\end{eqnarray}
As in the case $m>1$, with the help of   Lemma \ref{lem4.4}     we can show   that 
  the  quantity \eref{e5.11ii} has the  same $L^p$-limit as   
\begin{eqnarray}\label{e5.10}
 \frac12  n^{2H}  \Lambda_{ \eta(t) }   \sum_{e=2}^{5}  \int_{0}^{ \eta(t)  }
    \Gamma_{\eta(s)}
      dE_{e}(s) 
      &=& 
       \frac12  n^{2H}  \Lambda_{  \eta(t)}   \sum_{e=2}^{5}  \sum_{k=0}^{\lfloor \frac{nt}{T} \rfloor -1}  
    \Gamma_{t_{k}}
       \left(E_{e}(t_{k+1}) - E_{e}(t_{k }) \right).
\end{eqnarray}

 In the following, we show that the quantity in \eref{e5.10} converges to the solution of equation \eref{5.7}. Take $t\in \Pi$.
By \eref{eq4.10} we have
\begin{eqnarray*}
\sum_{k=0}^{  \frac{nt}{T}   -1}  
    \Gamma_{t_{k}}
       \left(E_{2}(t_{k+1}) - E_{2}(t_{k }) \right)
  &=&
  \frac14 \sum_{i =1}^d   \sum_{k=0}^{  \frac{nt}{T} -1}     \Gamma_{t_{k}}
\partial_{i} V   (X^n_{ t_k   })
    \int_{t_{k}}^{t_{k+1}} \partial V^{i}  (X^n_{v })      \left[ V(X^{n}_{t_{k+1} } ) + V(X^{n}_{t_{k}})  \right]  d B_{v}
    \\
    &&
    \quad\quad\quad\quad
\times \int_{t_k}^{t_{k+1} }  \int^s_{t_k}  dB_{u}
  dB_s.
\end{eqnarray*}
Take
\begin{eqnarray*}
 \wt{E}_{2}(t) &= &
   \frac12 \sum_{i =1}^d   \sum_{k=0}^{   \frac{nt}{T} -1}  \Gamma_{t_{k}}
(\partial_{i} V  \partial V^{i} V)  (X_{ t_k   })
    \int_{t_{k}}^{t_{k+1}}      d B_{v}
\int_{t_k}^{t_{k+1} }  \int^s_{t_k}  dB_{u}
  dB_s
  \\
 &= &
   \frac14 \sum_{i =1}^d   \sum_{k=0}^{   \frac{nt}{T} -1}    \Gamma_{t_{k}}
(\partial_{i} V  \partial V^{i} V)  (X_{ t_k   })
  ( B_{t_{k},t_{k+1}})^{3} .
  \end{eqnarray*}
It is easy to show that 
\begin{eqnarray}\label{eq4.33}
 n^{2H} \left(  \sum_{k=0}^{  \frac{nt}{T}   -1}  
    \Gamma_{t_{k}}
       \left(E_{2}(t_{k+1}) - E_{2}(t_{k }) \right) -  \wt{E}_{2}(t) \right) \rightarrow 0 \quad \text{ in }  L^{p} \quad  \text{ as } n \rightarrow \infty
\end{eqnarray}
  for $t\in \Pi$. 
Similarly, we  take
     \begin{eqnarray*}
\wt{E}_{3}(t) &= & 
 - \frac14   \sum_{i=1}^d  
 \sum_{k=0}^{   {nt}/{T}-1  }
    \Gamma_{t_{k}} \left(
  \partial (\partial_{i} V   V^{i}) V  
   + V^{i}   \partial( \partial_{i} V  )    V 
       \right)  ( X_{ t_{k}})  
 (  B_{t_{k}, t_{k+1}})^{3}\,,
\\
\wt{E}_{4} (t) &= & \wt{E}_{5} (t)=
\frac16 \sum_{i',i =1}^d \sum_{k=0}^{  nt/T -1  }    \Gamma_{t_{k}}(   V^{i'}  V^{i }  \partial_{i }\partial_{i'} V   ) (X_{t_{k}})   
( B_{t_{k},t_{k+1}})^{3}\,,
 \end{eqnarray*}
then   one can  show  that 
 \begin{eqnarray}\label{eq4.34}
n^{2H} \left( \sum_{k=0}^{  \frac{nt}{T}   -1}  
    \Gamma_{t_{k}}
       \left( E_{e}(t_{k+1}) - E_{e}(t_{k }) \right) - \wt{E}_{e}(t) \right) &\rightarrow& 0 
\end{eqnarray}
   { in } $L^{p}$          for    $e=3,\,4,\,5$.
  In summary from  \eref{eq4.33} and \eref{eq4.34}, we obtain
  \begin{eqnarray}\label{eq4.34i}
n^{2H} \sum_{e=2}^{5}  
 \sum_{k=0}^{  \frac{nt}{T}   -1}  
    \Gamma_{t_{k}}
       \left( E_{e}(t_{k+1}) - E_{e}(t_{k }) \right)   - n^{2H} \sum_{e=2}^{5} \wt{E}_{e}(t)  & \rightarrow& 0 
\end{eqnarray} 
    in $L^{p}$       for $t\in \Pi$.
 
 It is easy to see that
 \begin{eqnarray}\label{eq4.34ii}
\sum_{e=2}^{5}  \wt{E}_{e}(t) &= & -\frac 16 
 \sum_{i',i =1}^d \sum_{k=0}^{  nt/T -1  } \Gamma_{t_{k}} (   V^{i'}  V^{i }  \partial_{i }\partial_{i'} V   ) (X_{t_{k}})   
( B_{t_{k},t_{k+1}})^{3}
\end{eqnarray}
 for $t\in \Pi$. 
 Take $ f_{t} = \Gamma_{t} (   V^{i'}  V^{i }  \partial_{i }\partial_{i'} V   ) (X_{t })$ and $\zeta_{k,n} =  (B_{t_{k+1}} - B_{t_{k}})^{3}  $, then by  applying Proposition \ref{prop7.1} to  \eref{eq4.34ii}
 and taking into account Lemma \ref{lem6.1} (ii) 
    we obtain that
\begin{eqnarray*}
\frac12 \Lambda_{\eta(t)} \left(  n^{2H}  \sum_{e=2}^{5}   \wt{E}_{e}(\eta(t))  \right) \rightarrow   
   \bar{U}_{t} 
\end{eqnarray*}
  in $L^p$ for $t\in [0,T]$, where 
    \begin{eqnarray*}
 \bar{U}_{t} &= & 
  -  \frac{T^{2H}}{4}   \sum_{i',i =1}^d \Lambda_{t}  \int_{0}^{t} \Gamma_{s}(   V^{i'}  V^{i }  \partial_{i }\partial_{i'} V   ) (X_{t }) dB_{s}\,.
   \end{eqnarray*}
  Thanks to \eref{eq4.34i}, this convergence implies that
\begin{eqnarray*}
       \frac12  n^{2H}  \Lambda_{  \eta(t)}   \sum_{e=2}^{5}  \sum_{k=0}^{\lfloor \frac{nt}{T} \rfloor -1}  
    \Gamma_{t_{k}}
       \left(E_{e}(t_{k+1}) - E_{e}(t_{k }) \right) &\rightarrow & 
 \bar{U}_{t}
     \end{eqnarray*}
 for $t\in [0,T]$. 
The convergence \eref{5.6} follows from the fact that the process $\bar{U}$ solves  equation~\eref{5.7}.
   \hfill $\Box$

\section{Appendix}

\subsection{Proof of \eref{e3.14}}\label{section 6.1}
  The proof will be done in seven  steps.

\noindent \textit{Step 1.}\quad 
In this step, we derive a decomposition for $d_{2}$.  First,   applying  the integration by parts   formula (\ref{ipf}), we obtain
\begin{eqnarray}\label{e6.1}
\mE\left[
{{}{Z}}_n(t)
  \wt{D}_{u'}
   {{}{Z}}_n(t)
  D_{s'}
  {{}{Z}}_n(t)
  \right]
 & = &  
\sum_{k=0}^{ \lfloor \frac{nt}{T} \rfloor }
  \int_0^T\int_{0}^{T}\int_0^T\int_0^T
  \left[
 D_{r'} \wt{D}_{u'}
   {{}{Z}}_n(t)
   \right]
  \left[\wt{D}_{v'} D_{s'}
  {{}{Z}}_n(t)
  \right]
\nonumber  \\
& &\quad\quad\quad  \times  \beta_{\frac kn }(r)   \gamma_{t_k, r } (v ) 
   \mu(dvdv')  \mu(drdr')
\nonumber   \\
   &=  & 
   \sum_{k, k_3,k_4=0}^{ \lfloor \frac{nt}{T} \rfloor }
  \int_0^T\int_{0}^{T}\int_0^T\int_0^T
    \beta_{\frac {k_3}n} (r' ) \gamma_{t_{k_3}, r'  } (u' )
  \beta_{\frac {k_4}n} ( {s}' ) \gamma_{t_{k_4} ,  {s}'  } (v' )
  \nonumber
  \\
 &&\quad\quad\quad  \times \beta_{\frac kn }(r) \gamma_{t_k, r } (v )  
   \mu(dvdv')  \mu(drdr')
   ,
\end{eqnarray}
where the second equation follows from the fact that
   \begin{eqnarray*}
   \wt{D}_v D_r {Z}_n(t) &=&  \sum_{k=0}^{ \lfloor \frac{nt}{T} \rfloor  }
   \beta_{\frac kn} (r) \gamma_{t_k, r} (v), \quad t\in [0,T].
\end{eqnarray*}
Substituting the expression \eref{e6.1}  into \eref{eqn 3.8} we obtain
 \begin{eqnarray*}
  d_{2}&=&
  6 \sum_{k_1, k_2, k_3,k_4 =0}^{ \lfloor \frac{nt}{T} \rfloor }
  \int_{t_{k_4} }^{t_{k_4+1}}
  \int_{t_{k_1} }^{t_{k_1+1}}
  \int_0^T
  \int_0^T
    \int_{t_{k_3}}^{t_{k_3+1 } }
    \int_{t_{k_2} }^{t_{{k_2} +1}}
    \int_0^T
    \int_0^T
      \gamma_{t_{k_3}, r'  } (u' )
    \gamma_{t_{k_4} ,  {s}'  } (v' )
   \nonumber
    \\
    &&
      \quad\quad\quad \quad \quad \quad \quad \quad
   \times
     \gamma_{t_{k_2}, r } (v )
     \gamma_{t_{k_1}, s} (u)
    \mu(dvdv')  \mu(drdr')
   \mu(dudu')  \mu(dsds').
   \end{eqnarray*}
   By  changing the variables from $(v,v',r,r',u,u',s,s')$ to  $\frac{T}{n}(v,v',r,r',u,u',s,s')$   and   exchanging of the orders of   integrals associated with $\mu(dudu')$ and $\mu(drdr')$  we obtain
   \begin{eqnarray*}%\label{e 3.8}
 d_{2} &=&
  6  \left( \frac T n \right)^{8H}
  \sum_{k_1, k_2, k_3,k_4 =0}^{ \lfloor \frac{nt}{T} \rfloor }
 c(k_{1},k_{2},k_{3},k_{4}) ,
   \end{eqnarray*}
   where 
    \begin{eqnarray}\label{e6.2ii}
  c(k_{1},k_{2},k_{3},k_{4}) &=&
  \int_{ {k_4} }^{ {k_4+1}}
  \int_{ {k_1} }^{ {k_1+1}}
    \int_{ {k_3}}^{ {k_3+1 } }
    \int_{ {k_2} }^{ {{k_2} +1}}
   \int_0^n
  \int_0^n \int_0^n
    \int_0^n
      \varphi_{ {k_3}, r'  } (u' )
    \varphi_{ {k_4} ,  {s}'  } (v' )
   \nonumber
    \\
  &  &
     \quad \quad \quad 
   \times
    \varphi_{ {k_2}, r } (v )
     \varphi_{ {k_1}, s} (u)
     \mu(dvdv')  \mu(dudu') \mu(drdr')
   \mu(dsds'),
   \end{eqnarray}
and   recall that
     \begin{eqnarray}\label{e6.3ii}
  \varphi_{ k, s} (u)
      = \varphi^{0}_{ k, s} (u) - \varphi^{1}_{ k, s} (u), \quad 
   \varphi^{0}_{ k, s} (u) = \mathbf{1}_{ [k, s] }(u)  \quad  \text{ and}   \quad  \varphi^{1}_{ k, s} (u)= \mathbf{1}_{[s, {k+1}] }(u)   \,,
\end{eqnarray} 
where  $\mathbf{1}_{[a,b]}$ denotes  the indicator function of the interval $[a,b]$.

Now we denote
 \begin{eqnarray*}
I&:=&\left\{k_{1},k_{2},k_{3},k_{4} =0,1,\dots, \lfloor \frac{nt}{T} \rfloor \right\}.
\end{eqnarray*}
Take $i,j = 1,2,3,4 $, and denote by $I_{ij}$ the set of $(k_{1},k_{2},k_{3},k_{4})$ in $I$ such that $ |k_{i} - k_{j}|>2 $, that is,  $I_{ij} = \{(k_{1},k_{2},k_{3},k_{4}) \in I: |k_{i} - k_{j}|>2 \}$. Denote by $I_{ij}^{c}$ the complement of $I_{ij}$.  We   can  decompose $I$ as follows. 
\begin{eqnarray*}
I&=&   \bigcup_{l=1}^{8} M_{l} \,,
\end{eqnarray*}
 where
 \begin{eqnarray*}
M_{1} &=&  I_{42} \cap I_{41} \cap I_{31} \cap I_{32} ;
\\[3mm]
M_{2} &=&  \left( I_{42}^{c} \cap I_{41} \cap I_{31} \cap I_{32}\right) \bigcup  \left( I_{42} \cap I_{41} \cap I_{31}^{c} \cap I_{32}\right)
\\[1mm]
&:=& M_{21}+M_{22}
 ;
 \\[3mm]
M_{3} &=& \left( I_{42}  \cap I_{41}^{c} \cap I_{31} \cap I_{32}\right) \bigcup \left( I_{42}  \cap I_{41} \cap I_{31} \cap I_{32}^{c} \right);
 \\[3mm]
 M_{4} &=&  \left( I_{42}^{c}  \cap I_{41}^{c} \cap I_{31} \cap I_{32}\right) \bigcup  \left( I_{42}  \cap I_{41}^{c} \cap I_{31}^{c} \cap I_{32}\right)
 \\
 && \bigcup  \left( I_{42}^{c}  \cap I_{41}  \cap I_{31} \cap I_{32}^{c} \right) \bigcup  \left( I_{42}  \cap I_{41}  \cap I_{31}^{c} \cap I_{32}^{c } \right) 
 \\[1mm]
 &:=&   M_{41} \cup M_{42} \cup M_{43} \cup M_{44}   ;
 \\[3mm]
 M_{5} &=&  I_{42}  \cap I_{41}^{c} \cap I_{31} \cap I_{32}^{c}   ;
 \\[3mm]
 M_{6} &=&    I_{42}^{c}  \cap I_{41}  \cap I_{31}^{c} \cap I_{32}   ;
\\[3mm]
M_{7} &=&  \left( I_{42}^{c}  \cap I_{41}^{c} \cap I_{31}^{c} \cap I_{32}\right) \bigcup  \left( I_{42}^{c}  \cap I_{41}^{c} \cap I_{31}  \cap I_{32}^{c} \right)  
\\ 
&& \bigcup  \left( I_{42}^{c}  \cap I_{41}  \cap I_{31}^{c} \cap I_{32}^{c}\right) \bigcup  \left( I_{42}   \cap I_{41}^{c} \cap I_{31}^{c} \cap I_{32}^{c} \right)  ;
\\[3mm]
M_{8} &=&    I_{42}^{c}  \cap I_{41}^{c} \cap I_{31}^{c} \cap I_{32}^{c}   .
\end{eqnarray*}
  
  For any subset $M$ of $I$, we denote
\begin{eqnarray*}
 d_{2 } ( M)&:=&
  6  \left( \frac T n \right)^{8H}
  \sum_{(k_1, k_2, k_3,k_4) \in M } 
 c(k_{1},k_{2},k_{3},k_{4}) .
\end{eqnarray*}
  It is clear that 
  \begin{eqnarray*}
d_{2} &=& \sum_{l=1}^{8} d_{2}(M_{l}).
\end{eqnarray*}
Thus to show \eref{e3.14} it suffices to   show that $n^{8H-2} d_{2 } (M_{l}) \rightarrow 0$ as $n\rightarrow \infty$ for each $l=1,\dots, 8$.

\medskip

\noindent \textit{Step 2.}\quad  In this step, we   show the convergence of  $n^{8H-2} d_{2}(M_{7})$ and $n^{8H-2} d_{2}(M_{8})$. Since 
 \begin{eqnarray}\label{e6.3i}
| \varphi_{k,s} (u)| \leq \mathbf{1}_{[k,k+1]}(u),
\end{eqnarray}
 we have 
 \begin{eqnarray*}
|c ( k_{1},k_{2},k_{3},k_{4} )| &\leq & 1.
\end{eqnarray*}
Applying this inequality   to $d_{2}(M_{7})$ we obtain 
\begin{eqnarray*}
|d_{2}(M_{7})| &\leq &   6  \left( \frac T n \right)^{8H} \sum_{( k_{1},k_{2},k_{3},k_{4}) \in M_{7}} 1
 .
\end{eqnarray*}
 Note that  
\begin{eqnarray*}
M_{7} \subset \{|k_{i}-k_{j}|\leq 6  \text{ for }i,j=1,2,3,4 \},
\end{eqnarray*}
so the number of elements in $M_{7}$ is less than $6n$. 
 This implies that
  \begin{eqnarray*}
|d_{2}(M_{7})| 
&\leq &  36 n  \left( \frac T n \right)^{8H}.
\end{eqnarray*} 
 It follows from this estimate that $n^{8H-2} d_{2}(M_{7}) \rightarrow 0$ as $n\rightarrow 0$.  Note that $M_{8} \subset \{|k_{i}-k_{j}|\leq 4 \text{ for }i,j=1,2,3,4 \}$. So in the same way, we   can    show that $n^{8H-2} d_{2}(M_{8}) \rightarrow 0$.

\medskip

\noindent \textit{Step 3.}\quad In this step, we consider $d_{2}(M_{5})$ and $d_{2}(M_{6})$. Take $(k_{1},k_{2},k_{3},k_{4}  ) \in M_{5}$,  we have $|k_{2}-k_{4}|>2$ and $|k_{1}-k_{3}|>2$. By the mean value theorem and with the help of \eref{e6.3i}, it is easy to see that 
\begin{eqnarray}\label{e6.4}
|c( k_{1},k_{2},k_{3},k_{4})| &\leq& K |k_{2}-k_{4}|^{2H-2}|k_{1}-k_{3}|^{2H-2}.
\end{eqnarray}
Applying \eref{e6.4} to $d_{2}(M_{5})$ we obtain 
\begin{eqnarray}\label{e6.5i}
|d_{2}(M_{5})|&\leq &  K  \left( \frac T n \right)^{8H}
  \sum_{k_1, k_2, k_3,k_4 \in M_{5}}   |k_{2}-k_{4}|^{2H-2}|k_{1}-k_{3}|^{2H-2}
 .
\end{eqnarray}
Note that for $(k_1, k_2, k_3,k_4) \in M_{5}$ we have $ |k_{1} - k_{4}|\leq 2 $ and $ |k_{2} - k_{3}|\leq 2 $, so   
\begin{eqnarray*}
|k_{2} -k_{4}| &\leq& |k_{2}-k_{3}| + |k_{3} - k_{1}| + |k_{1} - k_{4}|
\\
&\leq & 3|k_{3} - k_{1}|.
\end{eqnarray*}
Applying this inequality to the right-hand side of \eref{e6.5i} yields 
\begin{eqnarray*}
| d_{2}(M_{5})| 
  &\leq&  K  \left( \frac T n \right)^{8H}
  \sum_{(k_1, k_2, k_3,k_4) \in M_{5}}   |k_{2}-k_{4}|^{2H-2}|k_{4}-k_{2}|^{2H-2}
  \\
  &\leq &  K  \left( \frac T n \right)^{8H}
  \sum_{ k_2,  k_4 :|k_{2}-k_{4}|>2}   |k_{2}-k_{4}|^{4H-4}.
\end{eqnarray*}
By taking $p=k_{2}-k_{4}$, we obtain
\begin{eqnarray*}
|d_{2}(M_{5})| &\leq&   K  \left( \frac T n \right)^{8H} \sum_{k_{2}  =0}^{n}   \sum_{n\geq |p|>2}  |p|^{4H-4}
\\
&\leq&   Kn  \left( \frac T n \right)^{8H} (n^{4H-3}\vee 1 ).
\end{eqnarray*}
It follows   from the above estimate that $n^{8H-2}d_{2}(M_{5})  $  converges to zero as $n$ tends to infinity. The proof for the convergence $n^{8H-2}d_{2}(M_{6}) \rightarrow 0$ is similar. Instead of  \eref{e6.4}, we  have the estimate
\begin{eqnarray*}
|c( k_{1},k_{2},k_{3},k_{4})| &\leq& K |k_{1}-k_{4}|^{2H-2}|k_{2}-k_{3}|^{2H-2}
\end{eqnarray*}
 for   $(k_{1},k_{2},k_{3},k_{4}) \in M_{6}$. 

\medskip

\noindent \textit{Step 4.}\quad In this step, we derive a new expression for $c( k_{1},k_{2},k_{3},k_{4})$.
Recall that $ \varphi_{ k_{4}, s'} (v')
      = \varphi^{0}_{ k_{4}, s'} (v') - \varphi^{1}_{ k_{4}, s'} (v')$.   Substituting this identity into \eref{e6.2ii} we obtain  
   \begin{eqnarray}\label{e6.3}
   c(k_{1},k_{2},k_{3},k_{4})&= & c_{0}-c_{1},
  \end{eqnarray}
  where
   \begin{eqnarray*} 
 c_{i} &=&   \int_{ {k_4} }^{ {k_4+1}}
  \int_{ {k_1} }^{ {k_1+1}}
    \int_{ {k_3}}^{ {k_3+1 } }
    \int_{ {k_2} }^{ {{k_2} +1}}
   \int_0^n
  \int_0^n \int_0^n
    \int_0^n
      \varphi_{ {k_3}, r'  } (u' )
    \varphi^{i}_{ {k_4} ,  {s}'  } (v' )
    \nonumber
    \\
    &&
    \quad \quad \quad \quad \quad \quad
   \times
    \varphi_{ {k_2}, r } (v )
     \varphi_{ {k_1}, s} (u)
    \mu(dvdv')  \mu(dudu') \mu(drdr')
   \mu(dsds')\,.
   \end{eqnarray*}
  By   exchanging  the orders of  the  integrals associated with $v'$ and $ s' $ in $c_{1}$,
we obtain
\begin{eqnarray*}
c_{1} &= &
\int_{ {k_4} }^{ {k_4+1}}
  \int_{ {k_1} }^{ {k_1+1}}
     \int_{ {k_3}}^{ {k_3+1 } }
    \int_{ {k_2} }^{ {{k_2} +1}}
   \int_{0}^{n} \int_{0}^{n}    \int_{  {k_4 }  }^{ v'}
      \int_{  0  }^{n }
        |v - v'|^{2H-2}  | {s} -  {s}'|^{2H-2} 
          \varphi_{ {k_3}, r'  } (u' ) 
 \\
  &&
   \quad\quad\quad\quad\quad\quad \times
   \varphi_{ {k_2}, r } (v )    \varphi_{ {k_1}, s} (u)
        dv ds'   \mu(dudu')  
  \mu( drdr' )
  ds dv' ,
\end{eqnarray*}
which, by switching the notations $s'$ and $v'$, is equal to
\begin{eqnarray*}
&
  \int_{ {k_4} }^{ {k_4+1}}
  \int_{ {k_1} }^{ {k_1+1}}
     \int_{ {k_3}}^{ {k_3+1 } }
    \int_{ {k_2} }^{ {{k_2} +1}}
  \int_{0}^{n} \int_{0}^{n} 
  \int_{  {k_4 }  }^{ s'}
      \int_{  0  }^{ n }
   |v - s' |^{2H-2}   | {s} -  v'|^{2H-2}     
       \varphi_{ {k_3}, r'  } (u' ) 
    \\
    & \quad\quad\quad\quad\quad \times      \varphi_{ {k_2}, r } (v )     \varphi_{ {k_1}, s} (u)
       dv dv'     \mu(dudu')  
         \mu( drdr' )
  ds ds' .
\end{eqnarray*}
Substituting the above expression of $c_{1}$ into \eref{e6.3} we obtain
\begin{eqnarray}
 c(k_{1},k_{2},k_{3},k_{4})  
 &=&
\int_{ {k_4} }^{ {k_4+1}}
  \int_{ {k_1} }^{ {k_1+1}}
   \int_{ {k_3}}^{ {k_3+1 } }
    \int_{ {k_2} }^{ {{k_2} +1}}
 \int_{0}^{n} \int_{0}^{n}     \int_{  {k_4 }  }^{ s'}
      \int_{  0  }^{ n }
    \phi (s,s',v , v' )   \varphi_{ {k_3}, r'  } (u' )        \nonumber
       \\
  &     &    \quad\quad\quad\quad
        \times
    \varphi_{ {k_2}, r } (v ) 
  \varphi_{ {k_1}, s} (u)
       dv dv'    \mu(dudu')  
  \mu( drdr' ) ds ds'
 ,
  \label{eqn 3.10}
  \end{eqnarray}
where  we denote
\begin{eqnarray*}
\phi (s,s',v , v' ) &=&
    | v -  v'|^{2H-2} |s - s' |^{2H-2}
       -  
     |v - s' |^{2H-2}    | {s} -  v'|^{2H-2} .
\end{eqnarray*}
 
 \medskip

\noindent \textit{Step 5.}\quad
We turn to  $d_{2}(M_{4})$.  It is easy to show    that 
\begin{eqnarray}\label{e6.5}
d_{2}  \left(M_{4i}\right) &=& d_{2}   \left( M_{4j}\right), \quad i,j =1,2,3,4.
\end{eqnarray}
 As an example, we show that $ d_{2}(M_{41}) =d_{2}(M_{44}) $. The other identities in \eref{e6.5}   can be shown similarly.  First, by exchanging the orders of integrals associated with $\mu(drdr')$ and $\mu(dsds')$ and integrals associated with $\mu(dvdv')$ and $\mu(dudu')$,  we obtain
 \begin{eqnarray*}
c(k_{1},k_{2},k_{3},k_{4}) &=&
    \int_{ {k_3}}^{ {k_3+1 } }
    \int_{ {k_2} }^{ {{k_2} +1}}  \int_{ {k_4} }^{ {k_4+1}}
  \int_{ {k_1} }^{ {k_1+1}}
   \int_0^n
  \int_0^n \int_0^n
    \int_0^n
      \varphi_{ {k_3}, r'  } (u' )
    \varphi_{ {k_4} ,  {s}'  } (v' )
   \nonumber
    \\
  &  &
     \quad \quad \quad 
   \times
    \varphi_{ {k_2}, r } (v )
     \varphi_{ {k_1}, s} (u)
   \mu(dudu')   \mu(dvdv') 
   \mu(dsds')  \mu(drdr').
\end{eqnarray*}
  Replacing $(v,v',u,u',r,r',s,s') $ by $(u,u',v,v',s,s',r,r')$ in the above expression, we obtain
   \begin{eqnarray*}
c(k_{1},k_{2},k_{3},k_{4}) &=&
    \int_{ {k_3}}^{ {k_3+1 } }
    \int_{ {k_2} }^{ {{k_2} +1}}  \int_{ {k_4} }^{ {k_4+1}}
  \int_{ {k_1} }^{ {k_1+1}}
   \int_0^n
  \int_0^n \int_0^n
    \int_0^n
      \varphi_{ {k_3}, s'  } (v' )
    \varphi_{ {k_4} ,  {r}'  } (u' )
   \nonumber
    \\
  &  &
     \quad \quad \quad 
   \times
    \varphi_{ {k_2}, s } (u )
     \varphi_{ {k_1}, r} (v)
   \mu(dvdv')   \mu(dudu') 
   \mu(drdr')  \mu(dsds')
   \\
   &=& c(k_{2},k_{1},k_{4},k_{3}).
\end{eqnarray*}
 So  we have
\begin{eqnarray*}
d_{2}(M_{41}) &=&   6  \left( \frac T n \right)^{8H}
  \sum_{(k_1, k_2, k_3,k_4 ) \in M_{41}} 
 c(k_{2},k_{1},k_{4},k_{3})
 = d_{2}(M_{44}),
\end{eqnarray*}
where the second identity follows by replacing $ (k_1, k_2, k_3,k_4 ) $ by $(k_{2},k_{1},k_{4},k_{3})$. 

The identities in \eref{e6.5} imply that to show the convergence $n^{8H-2}d_{2}(M_{4})\rightarrow 0$,  it suffices to show that  
$n^{8H-2}d_{2}(M_{44})\rightarrow 0$ as $n\rightarrow \infty$.

Take $(k_{1},k_{2},k_{3},k_{4}) \in M_{44}$. Then we have $|k_{1}-k_{4}|>2$ and $ |k_{2}-k_{4}|>2 $. This allows us to apply the mean value theorem to $\phi$ to obtain the estimate
 \begin{eqnarray}\label{e 3.19}
  \left|
    \phi(s,s', v, v')
     \right|
 &  \leq&
   K (|k_4-k_1|^{2H-3}|k_4-k_2|^{2H-2} + |k_4-k_1|^{2H-2 } |k_4-k_2|^{2H-3 } ).
    \end{eqnarray}
Applying \eref{e 3.19} to \eref{eqn 3.10} and taking into account \eref{e6.3i} we obtain
\begin{eqnarray}\label{e6.10}
 \left|
   c(k_1, k_2, k_3,k_4 )
     \right|
 &  \leq&
   K (|k_4-k_1|^{2H-3}|k_4-k_2|^{2H-2} + |k_4-k_1|^{2H-2 } |k_4-k_2|^{2H-3 } ).
\end{eqnarray}
Since $|k_{1}-k_{2}| \leq |k_{1}-k_{3}| + |k_{3}-k_{2}|  \leq 4 $, we have $|k_{4}-k_{1}| \leq 3 |k_{4}-k_{2}|  $.
This applied to \eref{e6.10} yields 
\begin{eqnarray*}
 \left|
   c(k_1, k_2, k_3,k_4 )
     \right|
 &  \leq&
K|k_{4}-k_{1}|^{4H-5},
\end{eqnarray*}
and thus
\begin{eqnarray*}
|d_{2}(M_{44}) |&\leq&   6  \left( \frac T n \right)^{8H} \sum_{( k_{1},k_{2},k_{3},k_{4}) \in M_{44}}   | c
(k_1, k_2, k_3,k_4 )|
\\
&\leq&    6  \left( \frac T n \right)^{8H} \sum_{k_{1},k_{4}: |k_{1}-k_{4}|>2}     K |k_4-k_1|^{4H-5} .
\end{eqnarray*}
By taking $p=k_{1}-k_{4}$, we obtain
\begin{eqnarray*}
|d_{2}(M_{44}) |&\leq&   6  \left( \frac T n \right)^{8H} \sum_{k_{1}  =0}^{n}   \sum_{n\geq |p|>2}  p^{4H-5}
\leq  12n  \left( \frac T n \right)^{8H} \sum_{p=3}^{\infty} p^{4H-5} ,
\end{eqnarray*}
which implies  that $n^{8H-2} d_{2}(M_{44}) \rightarrow 0$ as $n\rightarrow \infty$. 

  \medskip

\noindent \textit{Step 6.}\quad   In this step, we consider   $ d_{2}(M_{2})  $ and $d_{2}(M_{3}) $. 
As in Step 4, it is easy to show  that $d_{2}(M_{21}) = d_{2}(M_{22})$. So to show that $n^{8H-2}d_{2}(M_{2}) \rightarrow 0$  it suffices to show that $n^{8H-2}d_{2}(M_{22}) \rightarrow 0$ as $n\rightarrow \infty$.  

Take $(k_{1},k_{2},k_{3},k_{4}) \in M_{22}$, we     have $|k_{1}-k_{4}|>2$, $ |k_{2}-k_{4}|>2 $, and so the inequality \eref{e 3.19}   holds. 
Applying \eref{e 3.19}  to $d_{2}(M_{22})$ and taking $p_{1}=k_{1}-k_{4}$ and $p_{2}=k_{4}-k_{2}$, we obtain   
\begin{eqnarray*}
 |d_{2}(M_{22})| &\leq & K  \left( \frac T n \right)^{8H}
  \sum_{(k_1, k_2, k_{3}, k_4) \in M_{22}}     (|k_4-k_1|^{2H-3}|k_4-k_2|^{2H-2} + |k_4-k_1|^{2H-2 } |k_4-k_2|^{2H-3 } ) 
  \\
  &\leq & K  \left( \frac T n \right)^{8H}
  \sum_{k_{4}=1}^{n}\,\,\sum_{n \geq  |p_{1}|,\, |p_{2}| \geq 2 }   ( |p_{1}|^{2H-3} |p_{2}|^{2H-2}+ |p_{1}|^{2H-2} |p_{2}|^{2H-3} ) .
\end{eqnarray*}
It is   easy to see from the above estimate that $n^{8H-2}|d_{2}(M_{22})|  \leq Kn^{2H-2}$, which converges to zero as $n$ tends to infinity. 
The proof for the convergence $n^{8H-2}d_{2}(M_{3}) \rightarrow 0$ follows the same lines.

   \medskip

\noindent \textit{Step 7.}\quad    
    It remains to show that $n^{8H-2}d_{2}(M_{1}) \rightarrow 0$ as $n\rightarrow \infty$. To do this,  we first derive a new   expression for $c(k_{1},k_{2},k_{3},k_{4})$.
  Recall that 
       \begin{eqnarray*}
 \varphi_{ {k_3}, r'} (u')
&=& 
  \varphi_{ {k_3}, r'}^{0} (u') -   \varphi_{ {k_3}, r'}^{1} (u').
\end{eqnarray*}
  Substituting this identity     into \eref{eqn 3.10},  we obtain
    \begin{eqnarray}\label{e6.4i}
c(k_{1},k_{2},k_{3},k_{4}) &=& \tilde{c}_{0}-\tilde{c}_{1} ,
\end{eqnarray}
    where
    \begin{eqnarray*}
\tilde{c}_{i} &=&
\int_{ {k_4} }^{ {k_4+1}}
  \int_{ {k_1} }^{ {k_1+1}}
   \int_{ {k_3}}^{ {k_3+1 } }
    \int_{ {k_2} }^{ {{k_2} +1}}
 \int_{0}^{n} \int_{0}^{n}     \int_{  {k_4 }  }^{ s'}
      \int_{  0  }^{ n }
    \phi (s,s',v , v' )   \varphi_{ {k_3}, r'  }^{i} (u' )        \nonumber
       \\
  &     &    \quad\quad\quad\quad
        \times
    \varphi_{ {k_2}, r } (v ) 
  \varphi_{ {k_1}, s}   (u)
       dv dv'    \mu(dudu')  
  \mu( drdr' ) ds ds'\,.
\end{eqnarray*}

As in Step 3, by exchanging the   order of the integrals associated with the variables $r'$ and $u'$, and then switching the notations $r'$ and $u'$, we obtain
\begin{eqnarray*}
\tilde{c}_{1}  &=&
  \int_{ {k_4} }^{ {k_4+1}}
  \int_{ {k_1} }^{ {k_1+1}} \int_{ {k_3}}^{ {k_3+1 } }
    \int_{ {k_2} }^{ {{k_2} +1}}
     \int_{  k_{3}  }^{r' } \int_{ 0 }^{n }
 \int_{  {k_4 }  }^{ s'}
      \int_{ 0  }^{n }
     \phi (s  ,s',v , v' )
     |u-r'|^{2H-2}  |r-u'|^{2H-2}
      \\
     & &\quad\quad\quad\quad 
      \times    \varphi_{ {k_2}, r  } (v )    \varphi_{ {k_1}, s  } (u)    dv dv'
  d ud {u}'
    drdr'  ds ds'
  \,.
  \end{eqnarray*}
Substituting the above expression of $\tilde{c}_{1}$ into \eref{e6.4i}, we obtain
\begin{eqnarray}\label{e 3.18}
c (k_{1},k_{2},k_{3},k_{4})&= &
\int_{ {k_4} }^{ {k_4+1}}
  \int_{ {k_1} }^{ {k_1+1}}
  \int_{ {k_3}}^{ {k_3+1 } }
    \int_{ {k_2} }^{ {{k_2} +1}}    \int_{  k_{3}  }^{r' } \int_{  0  }^{n }   
       \int_{  {k_4 }  }^{ s'}
      \int_{  0  }^{n }
  \phi(s,s', v, v')  \phi(r,r',u,u')
   \nonumber
  \\
  &&\quad\quad\quad\quad\quad
  \times
        \varphi _{ {k_2}, r } (v )     \varphi_{ {k_1}, s  } (u ) dv dv' d {u}d u'
  drdr'   ds ds'
 .
\end{eqnarray}
%where
 %\begin{eqnarray*}
%\wt{\phi} ( s,s',v,v',r,r',u,u') &=&   \phi(s,s', v, v')  \phi(r,r',u,u')
 %  .
%\end{eqnarray*}

Take $(k_{1},k_{2},k_{3},k_{4}) \in M_{1}$, then it is clear that the inequality \eref{e 3.19} holds true, and in the same way, we can show that
\begin{eqnarray}\label{e6.11i}
 | \phi(r,r',u,u') |&\leq & K( |k_{1}-k_{3}|^{2H-3}|k_{2}-k_{3}|^{2H-2}+ |k_{1}-k_{3}|^{2H-2}|k_{1}-k_{3}|^{2H-3}). 
\end{eqnarray}
Applying   inequalities \eref{e 3.19} and \eref{e6.11i}   to \eref{e 3.18} and taking $p_{1}=k_{3}-k_{1}$, $p_{2} = k_{2}-k_{3}$, $p_{3} = k_{4} - k_{2}$     we   obtain 
\begin{eqnarray}\label{e6.12i}
| d_{2}(M_{1}) | &\leq&   6  \left( \frac T n \right)^{8H}
  \sum_{(k_1, k_2, k_3,k_4 ) \in M_{1}} 
 |c(k_1, k_2, k_3,k_4)|
 \nonumber
 \\
 &\leq& Kn^{-8H}   \sum_{k_{1}=0}^{n} \sum_{(p_{1},p_{2},p_{3}) \in J}  ( |p_{1}|^{2H-3}|p_{2}|^{2H-2}+ |p_{1}|^{2H-2}|p_{2}|^{2H-3})
 \tilde{c} (p_{1},p_{2},p_{3})  
\nonumber  \\
 &=& Kn^{-8H}   \sum_{k_{1}=0}^{n} \sum_{p_{1},p_{2}: n \geq |p_{1}|,|p_{2}| >2}  ( |p_{1}|^{2H-3}|p_{2}|^{2H-2}+ |p_{1}|^{2H-2}|p_{2}|^{2H-3})
\nonumber
 \\
 && \quad\quad\quad\quad\quad\quad\quad\quad\quad\quad\quad\quad \times \sum_{ p_{3}  \in J(p_{1},p_{2})}    \tilde{c} (p_{1},p_{2},p_{3})
 ,
\end{eqnarray}
where 
\begin{eqnarray*}
\tilde{c} (p_{1},p_{2},p_{3}) &=& |p_{1}+p_{2}+p_{3}|^{2H-3}|p_{3}|^{2H-2} +|p_{1}+p_{2}+p_{3}|^{2H-2}|p_{3}|^{2H-3}  ,
\end{eqnarray*}
\begin{eqnarray*}
J = \{(p_{1},p_{2},p_{3}): n \geq |p_{1}|,|p_{2}| , |p_{3}|, |p_{1}+p_{2}+p_{3}|>2\},
\end{eqnarray*}
and
  \begin{eqnarray*}
J(p_{1},p_{2})= \{ p_{3}: n \geq |p_{3}|, |p_{1}+p_{2}+p_{3}|> 2 \}.
\end{eqnarray*}

We claim that $\sum_{p_{3} \in J(p_{1},p_{2})} \tilde{c} (p_{1},p_{2},p_{3}) $ is uniformly bounded in $(p_{1},p_{2})$ by a constant. 
 Take $p_{1},p_{2}$ such that $n \geq  |p_{1}|,|p_{2}|>2$.  Since when $|p_{1}+p_{2} +p_{3}|< |p_{3}|$  we have $|p_{1}+p_{2} +p_{3}|^{\al}>  |p_{3}|^{\al}$ for $\al=2H-2$ and $\al=2H-3$. So 
\begin{eqnarray}\label{e6.14i}
 \sum_{p_{3}\in J(p_{1},p_{2}):  |p_{1}+p_{2} +p_{3}| < |p_{3}|}  
 \tilde{c} (p_{1},p_{2},p_{3}) 
  &\leq&
  \sum_{p_{3}\in J(p_{1},p_{2}) } |p_{1}+p_{2}+p_{3}|^{4H-5}
  \nonumber
  \\
  &\leq&2 \sum_{p =3}^{\infty}p^{4H-5}.
\end{eqnarray}
Similarly, we have
\begin{eqnarray}\label{e6.15i}
 \sum_{p_{3}\in J(p_{1},p_{2}): |p_{1}+p_{2} +p_{3}|\geq |p_{3}|}   \tilde{c} (p_{1},p_{2},p_{3})
  &\leq&
  \sum_{p_{3}:n \geq |p_{3}| >2} |p_{3}|^{4H-5}
  \leq2\sum_{p=3}^{\infty} p^{4H-5}.
\end{eqnarray}
In summary of \eref{e6.14i} and \eref{e6.15i}, we have shown that 
\begin{eqnarray}\label{e6.17}
 \sum_{p_{3}\in J(p_{1},p_{2}) }   \tilde{c} (p_{1},p_{2},p_{3})
   &\leq& 4\sum_{p=2}^{\infty} p^{4H-5}.
\end{eqnarray}
 
 Applying   inequality \eref{e6.17} to \eref{e6.12i}, we obtain the estimate
 \begin{eqnarray*}
| d_{2}(M_{1}) | &\leq&K n^{-6H},
\end{eqnarray*}
which implies that $n^{8H-2}d_{2}(M_{1}) \rightarrow 0$ as $n\rightarrow \infty$. \hfill $\Box$

\subsection{Proof of \eref{e3.34}}\label{section 6.2}

By the integration by parts  formula (\ref{ipf1}), we obtain
\begin{eqnarray*}
\mE(z_{n}(t) B_{r}) &=&  \sum_{k=0}^{ \lfloor \frac{nt}{T} \rfloor  } \int_{0}^{r} \int_{0}^{T} \int_{0}^{T}      \beta_{\frac kn} (s)    \gamma_{t_k, s} (u)     d{u} \mu(d{s} ds').
\end{eqnarray*}
By changing the variables from $(u,s,s')$ to $\frac{T}{n} (u,s,s')$ in the above expression, we obtain
\begin{eqnarray}\label{e6.19}
\mE(z_{n}(t) B_{r}) &=&\left(\frac{T}{n}\right)^{2H+1}   \sum_{k=0}^{ \lfloor \frac{nt}{T} \rfloor  }  \int_{0}^{\frac{nr}{T}}  \int_{k}^{k+1}  \int_{0}^{n}        \varphi_{ k, s} (u)   |s-s'|^{2H-2}  d{u}  d{s} ds'\,,
\end{eqnarray}
where $\varphi_{k,s}(u)$ is defined  in \eref{e6.3ii}. Let us denote $I_{1}(k) = [k-2, k+2] \cap [0, \frac{nr}{T}]$ and $I_{2}(k) = [0, \frac{nr}{T}] \setminus I_{1}(k) $, and set
\begin{eqnarray*}
A_{i} &=&\left(\frac{T}{n}\right)^{2H+1}   \sum_{k=0}^{ \lfloor \frac{nt}{T} \rfloor  }  \int_{I_{i}(k)}   \int_{k}^{k+1}  \int_{0}^{n}        \varphi_{ k, s} (u)   |s-s'|^{2H-2}  d{u}  d{s} ds'\,.
\end{eqnarray*}
 Then it is easy to show that
 \begin{eqnarray*}
\mE(z_{n}(t) B_{r}) &=& A_{1}+A_{2}. 
\end{eqnarray*}
So to show that $n^{H+\frac12} \mE(z_{n}(t) B_{r})  \rightarrow 0$  it suffices to show that $n^{H+\frac12} A_{i} \rightarrow 0$ as $n\rightarrow \infty$ for $i=1,2$.

We can write
\begin{eqnarray*}
A_{1}&\leq&  \left(\frac{T}{n}\right)^{2H+1}   \sum_{k=0}^{ \lfloor \frac{nt}{T} \rfloor  }  \int_{k-2}^{k+2}   \int_{k}^{k+1}  \int_{k}^{k+1}        |s-s'|^{2H-2}  d{u}  d{s} ds'
\leq K n  \left(\frac{T}{n}\right)^{2H+1} ,
\end{eqnarray*}
so we have $n^{H+\frac12} A_{1} \rightarrow 0$ as $n\rightarrow \infty$.

Now we turn to $A_{2}$. By  exchanging  the orders of the integrals with respect to $u$ and $s$, we have 
\begin{eqnarray}\label{e6.20}
\int_{k}^{k+1}  \int_{0}^{n}        \mathbf{1}_{ [k, s]} (u)   |s-s'|^{2H-2}  d{u}  d{s}  &=& \int_{k}^{k+1}  \int_{k}^{s}          |u-s'|^{2H-2}  d{u}  d{s} .
\end{eqnarray}
 Substituting \eref{e6.20} into $A_{2}$ we obtain 
\begin{eqnarray}\label{e6.21}
A_{2} &=&\left(\frac{T}{n}\right)^{2H+1}   \sum_{k=0}^{ \lfloor \frac{nt}{T} \rfloor  }  \int_{I_{2}(k)}   \int_{k}^{k+1}  \int_{k}^{s}       ( |s-s'|^{2H-2} - |u-s'|^{2H-2} )   d{u}  d{s} ds'\,.
\end{eqnarray}
 Note that for $s' \in I_{2}(k)$   we have
\begin{eqnarray*}
\left|  \int_{k}^{k+1}  \int_{k}^{s}       ( |s-s'|^{2H-2} - |u-s'|^{2H-2} )   d{u}  d{s} \right| &\leq& K |k-s'|^{2H-3},
\end{eqnarray*}
so
\begin{eqnarray*}
|A_{2}| &\leq&  \left(\frac{T}{n}\right)^{2H+1}   \sum_{k=0}^{ \lfloor \frac{nt}{T} \rfloor  }  \int_{I_{2}(k)} K  |k-s'|^{2H-3} ds'
\\
 &=& K \left(\frac{T}{n}\right)^{2H+1}   \sum_{k=0}^{ \lfloor \frac{nt}{T} \rfloor  } (  2^{2H-2}-( n-k )^{2H-2} + 2^{2H-2}- k^{2H-2}  ) 
 \\
 &\leq & K \left(\frac{T}{n}\right)^{2H+1} \sum_{k=0}^{n}   2^{2H-2} = Kn \left(\frac{T}{n}\right)^{2H+1}  ,
\end{eqnarray*}
which implies that $n^{H+\frac12} A_{2} \rightarrow 0$ as $n \rightarrow \infty$.
This completes the proof. \hfill $\Box$
\subsection{Estimates of some triple integrals}
In this subsection, we provide   estimates for some triple integrals which have been used the main body of the paper.

 \begin{lemma}\label{lem6.1}
(i) For $t\in \Pi$,  we define
\begin{eqnarray}\label{eq6.1}
G_t  & =  & \sum_{k=0}^{  nt/T   -1 } \int_{t_{k}}^{t_{k+1}} \int_{t_{k}}^{s } \int_{t_{k}}^{u}   dB^1_{v} dB^2_{u} dB^3_{s },
\end{eqnarray}
 where    $B^{1}$, $B^{2}$, $B^{3}$   are either a fractional Brownian motion with Hurst parameter $H>\frac12$ or equal to the identity function. Take $p \geq 1$. Then we have
  \begin{eqnarray*}
 \|   {G}_t - G_s \|_p    \leq K n^{-2 H} |t-s|^{\frac12}
\end{eqnarray*}
for $s,t \in \Pi$.   

(ii) Let $B$ a one-dimensional fBm with Hurst parameter $H>\frac12$.  Take $p\geq1$ and   $t \in [0, T]$. We have the following convergence in $L^p$:
\begin{eqnarray*}
n^{2H} \sum_{k=0}^{\lfloor nt/T \rfloor-1} (B_{ t_{k+1} } - B_{t_{k}  })^{3}  &\rightarrow& 3 {T^{2H} }  B_t\,.
\end{eqnarray*}

  \end{lemma}
 \noindent \textit{Proof:}\quad 
  The result in (i)  follows from Proposition 5.10 in \cite{HLN2}. The convergence in (ii) follows immediately from results  in \cite{GS} or \cite{T}. 
 \hfill $\Box$
  
  \medskip
  
  We need the following technical lemma.
  
\begin{lemma}\label{lem5.2}
Let $f$ and $g$ be $\beta$-H\"older continuous stochastic processes   on $[0,T]$ and $ \mE[ \|f\|_{\beta}^{p}]+ \mE[ \|g\|_{\beta}^{p}] \leq K$ for all $\frac12 < \beta < H$ and $p\geq 1$, and let  $h^{n}$, $n\in \NN$ be   processes on $[0,T]$ such that 
\begin{eqnarray*}
\|h^{n}_{t } - h^{n}_{s}\|_{p}& \leq & K (t-s)^{\beta}, \quad s,t \in \Pi: s\leq t.
\end{eqnarray*}
Let $B^{1}$, $B^{2}$ and $B^{3}$ be  as in Lemma \ref{lem6.1}.   %Take $i,j=1$, $2$ or $3$, and 
Define the process 
\begin{eqnarray}\label{e6.23}
\wt{G}_{t}  &= & \sum_{k=0}^{  nt/T  -1 } h^{n}_{t_{k}}  \int_{t_{k}}^{t_{k+1}} \int_{t_{k}}^{s_{3}} \int_{t_{k}}^{s_{2}} f_{s_{i}} g_{s_{j}}   dB^1_{s_{1}} dB^2_{s_{2}} dB^3_{s_{3}}\,,\quad  t  \in \Pi\,. 
\end{eqnarray}
Then the following estimate holds true  for all $s,t \in \Pi$:
  \begin{eqnarray}\label{e6.2i}
 \|   \wt{G}_t - \wt{G}_s \|_p   & \leq  &K n^{-2 H} |t-s|^{\frac12}.
\end{eqnarray}
\end{lemma}
 \noindent \textit{Proof:}\quad 
  We decompose $\wt{G}$ as follows:
  \begin{eqnarray}\label{e6.2}
\wt{G}_{t}  & =  & \sum_{k=0}^{  nt/T  -1} f_{t_{k}} g_{t_{k}} h^{n}_{t_{k}} \int_{t_{k}}^{t_{k+1}} \int_{t_{k}}^{s_{3}} \int_{t_{k}}^{s_{2}}  dB^1_{s_{1}} dB^2_{s_{2}} dB^3_{s_{3}}
\nonumber
\\
  &&+ \sum_{k=0}^{  nt/T  -1} f_{t_{k}}   h^{n}_{t_{k}} \int_{t_{k}}^{t_{k+1}} \int_{t_{k}}^{s_{3}} \int_{t_{k}}^{s_{2}} \int_{t_{k}}^{s_{j}} dg_{s_{4}} dB^1_{s_{1}} dB^2_{s_{2}} dB^3_{s_{3}}
  \nonumber
\\
&& +   \sum_{k=0}^{  nt/T  -1} h^{n}_{t_{k}} \int_{t_{k}}^{t_{k+1}} \int_{t_{k}}^{s_{3}} \int_{t_{k}}^{s_{2}}  \int_{t_{k}}^{s_{i}}  g_{s_{j}}  d f_{s_{4}} dB^1_{s_{1}} dB^2_{s_{2}} dB^3_{s_{3}}\,.
\end{eqnarray}
Applying Proposition 5.10 in \cite{HLN2} to  the second and third terms on the right-hand side of \eref{e6.2}, and applying Lemma \ref{lem 2.4} to the first term and taking into account the estimate in Lemma \ref{lem6.1} (i), we obtain the inequality \eref{e6.2i}.
 \hfill $\Box$

 \subsection{Proof of Lemma \ref{lem4.4}}\label{section 6.3i} 

 By the definition of $J_{1}$, we have
\begin{eqnarray*}
  \int_{s}^{ t }
    \Gamma^{n }_{u}
      dJ_{1}(u)
      &= &  
  \sum_{i=1}^{d}     \int_{s}^{ t }
    \Gamma^{n }_{u} 
   \int_{\eta(u)}^{u} \partial_{i} V (X^{n}_{v}) d X^{n,i}_{v}    dB_u
   \\
   &= &
    \sum_{i=1}^{d}   \sum_{k=\lfloor \frac{ns}{T} \rfloor}^{\lfloor \frac{nt}{T} \rfloor}  \int_{t_{k}\vee s}^{ t_{k+1}\wedge t }
    \Gamma^{n }_{u} 
   \int_{\eta(u)}^{u} \partial_{i} V (X^{n}_{v}) d X^{n,i}_{v}    dB_u
\end{eqnarray*}
for $   s,t \in [0,T]$.
Applying the Minkovski inequality to the right-hand side of the above equation and taking into account   Lemma 8.2 in \cite{HLN2}  we obtain the estimate
\begin{eqnarray}\label{e 5.19}
\left\| \int_{s}^{ t }
    \Gamma^{n }_{u}
      dJ_{1}(u)
\right\|_{p} & \leq &
  \sum_{i=1}^{d}   \sum_{k=\lfloor \frac{ns}{T} \rfloor}^{\lfloor \frac{nt}{T} \rfloor} n^{-\beta} (  t_{k+1}\wedge t  -  t_{k}\vee s )^{\beta}
  \nonumber
\\
& \leq & K |t-s|^{\beta} n^{1-2\beta},
\end{eqnarray}
for $i=1$. In the same way we can show       that   estimate \eref{e 5.19} holds while $J_{1}$ is replaced by $J_{2}$. Applying these two estimates to 
  \begin{eqnarray*}
Y_t &= &
 \frac12  \Lambda^{n }_{  t}   \sum_{i=1}^{2}  \int_{0}^{ t }
    \Gamma^{n }_{s}
      dJ_{i}(s)\,
\end{eqnarray*}
 we obtain 
\begin{eqnarray}\label{eq4.25}
\|Y\|_{\beta, p} &\leq & K n^{1-2\beta}.
\end{eqnarray}

We  denote $\Phi := \Lambda - \Lambda^{n } $. Subtracting \eref{eqn 3.1i} from \eref{e 3.3} we obtain
\begin{eqnarray*}
\Phi^{i}_{i'}  (t)
&=
  &
  \sum_{j=0}^m \sum_{i'' =1}^d\int_0^t
 \left[
 \partial_{i''} V^i_j (  X_s  )\Lambda^{ i''}_{i'}({s})
 -
 V^i_{j  i''} (s) \Lambda^{n,i''}_{i'}( {s})
 \right]
 dB^{j }_s
 \\
 &=&
  \sum_{j=0}^m \sum_{i'' =1}^d\int_0^t
  \partial_{i''} V^i_j (  X_s  )
 \Phi^{i''}_{i'}({s})
 dB^{j }_s
  +
  \sum_{j=0}^m \sum_{i'' =1}^d\int_0^t
 \left[
 \partial_{i''} V^i_j (  X_s  )
 -
 V^i_{j  i''}({s})
 \right]
 \Lambda^{n,i''}_{i'} ({  s})
 dB^{j }_s \,.
   \end{eqnarray*}
By     the product rule it is easy to verify the following   identity,
   \begin{eqnarray}\label{eqn 3.15ii}
\Lambda (t) - \Lambda^{n} (t)
&=&
   \sum_{ i,i'=1}^d
  \sum_{j=0}^m     \Lambda (t)
   \int_0^t
\Gamma_{i'}(s) 
 \left[
 \partial_{i} V^{i'}_j (  X_s  )
 -
 V^{i'}_{j  i}({s})
 \right]
 \Lambda^{n,i} (s)
 dB^{j }_s \,.
   \end{eqnarray}
Denote 
 \begin{eqnarray*}
\wt{V}(X_{s}, X_{s}^{n}) &= &  
\int_{0}^{1}  \int_{0}^{1} \partial  \partial_{i'} V^{i''}_{j} \Big(\lambda X_{s} + (1-\lambda) (\theta X_{s} +(1-\theta) X^{n}_{s} ) \Big) (1-\theta) d\lambda d\theta   \,.
\end{eqnarray*}
It is easy to verify that 
 \begin{eqnarray*}
 \partial_{i'} V^{i''}_j (  X_s  )
 -
 V^{i''}_{j  i'}({s}) &= & \wt{V}(X_{s}, X^{n}_{s} ) Y_{s} \, .
\end{eqnarray*}
Then \eref{eqn 3.15ii} becomes
\begin{eqnarray*}
  \Lambda_{i}({t}) -\Lambda^{n}_{i}({t})     &= & 
\Lambda_{t} \int_0^t
 dg_s \cdot Y_{s} \,,
\end{eqnarray*}  
where  
\begin{eqnarray*}
 g_{t} &= &
  \sum_{ i',i''=1}^d
  \sum_{j=0}^m     \int_0^t
\Gamma_{i''}(s) 
 \wt{V}(X_{s}, X^{n}_{s} )  
 \Lambda^{n,i'}_{  i}(s)
 dB^{j }_s \,.
\end{eqnarray*}
Applying Lemma \ref{lem 2.3} and taking into account the estimate \eref{eq4.25}, we obtain 
 \begin{eqnarray*}
\left\| \Lambda_{t} \int_s^t
 dg_s \cdot Y_{s}\right\|_{p} &\leq & K n^{1-2\beta} (t-s)^{\beta},
\end{eqnarray*}
 which implies the estimate for $\| \Lambda-\Lambda \|_{\beta, p}$.
   The  estimate for the quantity  $\Gamma-\Gamma^{n}$  can be shown similarly.  
\hfill$\Box$

\subsection{Proof of \eref{e4.23}}\label{section 6.3}

It is clear that
\begin{eqnarray}\label{e6.11}
&&
\sum_{i=1}^{2 } \int_{0}^{t } \Gamma^{n}_{u} d J_{i}(u) - \sum_{e=2}^{5 }\int_{0}^{\eta(t)} \Gamma^{n}_{\eta(u)} d E_{e}(u)
\nonumber
\\
&&=
\sum_{i=1}^{2 } \int_{0}^{t } \Gamma^{n}_{u} d J_{i}(u) - \sum_{i=1}^{2 }\int_{0}^{\eta(t)} \Gamma^{n}_{\eta(u)} d J_{i}(u)
\nonumber
\\
&&= 
\sum_{i=1}^{2 } \int_{0}^{t } \int_{\eta(s)}^{s} d\Gamma^{n}_{u} d J_{i}(s) + \sum_{i=1}^{2 }\int_{\eta(t)}^{t} \Gamma^{n}_{\eta(u)} d J_{i}(u).
\end{eqnarray}
In the following, we estimate  the $L^p$-norms of the two terms on the right-hand side of \eref{e6.11}.

For $t \in [0,T]$, we define
\begin{eqnarray*}
I   (t) %=\,& \frac12 
   %\sum_{k=0}^{  \lfloor {nt}/{T} \rfloor -1}
  % \partial  V   V  (X^n_{t_k})
  % (B_{t_{k+1}\wedge t} - B_{t_{k}})^{2}
   %\\
   &= & \int_{0}^{t}  (\partial  V   V ) (X^n_{\eta(s)}) (B_{s} - B_{\eta(s)})   dB_{s} ,
      \end{eqnarray*}
     It is clear that $I(t_{k})=I_{1}(t_{k})=I_{2}(t_{k})$ for $k=0,1,\dots, n$. As in \eref{eq4.5}, we take the decomposition
     \begin{eqnarray*}
J_{1}(t) + J_{2}(t)  &= &
         \left( R_0 (t) - I  (t) \right)
      +\left( I  (t)-\wt{R}_{0} (t) \right)  +R_1(t) +\wt{R}_1(t) 
    \nonumber
      \\
     & := &
       E_2 (t) +E_3 (t) +E_4 (t) +E_5 (t) 
\end{eqnarray*}
     for $t \in [0,T]$, where $R_{0}, R_{1}, \wt{R}_{0}$ and $\wt{R}_{1}$ are defined as before. Note that the $E_{2}$ and $E_{3}$ defined here are extensions of those in \eref{eq4.5} from $\Pi$ to   $[0,T]$. 
     
     By applying Lemma 8.2 in \cite{HLN2}  to    \eref{e4.5ii} and \eref{e4.5i} we obtain
     \begin{eqnarray}\label{e6.12}
\|E_{e}\|_{[t_{k}, t_{k+1}], \beta} &\leq & Ke^{K\|B\|_{\beta}^{1/\beta}} n^{ -2\beta}, \quad e=4,5.
\end{eqnarray}
  Similarly, we can show that inequality \eref{e6.12}  also holds for $e=2,3$. Therefore, we obtain
\begin{eqnarray}\label{e6.13}
\|J_{1}+J_{2}\|_{[t_{k}, t_{k+1}],  \beta} &= &
\left\|\sum_{e=2}^{5}E_{e} \right\|_{[t_{k}, t_{k+1}],  \beta}
\nonumber
\\ 
 &\leq  & Ke^{K\|B\|_{\beta}^{1/\beta}}  n^{-2\beta}.
\end{eqnarray}
     By applying Lemma 8.2 in \cite{HLN2} and with the help of the estimate \eref{e6.13} we obtain
      \begin{eqnarray*}
\left\|
 \sum_{i=1}^{2} \int_{t'}^{t''} \int_{\eta(s)}^{s} d \Gamma^{n}_{u} dJ_{i}(s)
\right\|_{p} &\leq & Kn^{ -4\beta}
\end{eqnarray*}
for $t', t'' \in [t_{k},t_{k+1}]$.
     Therefore, we have
     \begin{eqnarray*}
\left\|
  \sum_{i=1}^{2} \int_{0}^{t} \int_{\eta(s)}^{s} d \Gamma^{n}_{u} dJ_{i}(s)
\right\|_{p} 
&\leq &
\sum_{k=0}^{\lfloor nt/T \rfloor}
\left\|
  \sum_{i=1}^{2} \int_{t_{k}}^{t_{k+1}\wedge t} \int_{\eta(s)}^{s} d \Gamma^{n}_{u} dJ_{i}(s)
\right\|_{p} 
\\
 & \leq & Kn^{ 1-4\beta}.
\end{eqnarray*}  
 On the other hand,  applying \eref{e6.13} to
     \begin{eqnarray*}
   \sum_{i=1}^{2} \int_{\eta(t)}^{t}  \Gamma^{n}_{\eta(s)} dJ_{i}(s)  
&=& 
   \sum_{i=1}^{2}    \Gamma^{n}_{\eta(t)} ( J_{i}(t)  - J_{i}(\eta(t)) \,,
\end{eqnarray*}
     we   obtain
     \begin{eqnarray*}
\left\|
  \sum_{i=1}^{2} \int_{\eta(t)}^{t}  \Gamma^{n}_{\eta(s)} dJ_{i}(s)  \right\|_{p}
 \leq  Kn^{-3\beta}\leq \, Kn^{1-4\beta}.
\end{eqnarray*}
     This  completes the proof. \hfill $\Box$

\end{document}